\input amstex
\documentstyle{amsppt}
\magnification=1200

\def\SetGrid#1#2{
\dimen1= #1      
\dimen2= #2      
\dimen0=\dimen1  
\dimen3=\dimen1  
\dimen4=\dimen2  
\dimen6=\dimen2  %
}

\def\SetOrigin #1 #2 {      
\dimen7= \dimen1        
\multiply\dimen7 by #1
\dimen8= \dimen2             
\multiply\dimen8 by -#2
}

\def\LineWidth#1{\dimen5=#1}

\def\li #1 #2 #3 {\vrule height #1  width #2  depth #3 }

\def\hl #1 #2 #3 {
\ifnum #3 > 0 
{
\multiply\dimen1 by #1
\advance\dimen1 by \dimen7
\multiply\dimen2 by -#2
\advance\dimen2 by \dimen8
\multiply\dimen3 by #3
\rlap{ \kern \dimen1
\raise \dimen2 
\hbox{ \li {\dimen5} {\dimen3} 0pt }}
}
\else
{
\multiply\dimen1 by #1
\advance\dimen1 by \dimen7
\multiply\dimen2 by -#2
\advance\dimen2 by \dimen8
\multiply\dimen3 by -#3
\advance\dimen1 by -\dimen3
\rlap{ \kern \dimen1
\raise \dimen2 
\hbox{ \li {\dimen5} {\dimen3} 0pt }}
}
\fi}


\def\vl #1 #2 #3 {
\ifnum #3 > 0 
{
\multiply\dimen1 by #1
\multiply\dimen2 by -#2
\multiply\dimen4 by #3
\advance\dimen1 by \dimen7
\advance\dimen2 by \dimen8
\advance\dimen2 by -\dimen4
\rlap{ \kern \dimen1 \raise \dimen2 
\hbox{ \li {\dimen4} {\dimen5} 0pt }}
}
\else
{
\multiply\dimen1 by #1
\multiply\dimen2 by -#2
\advance\dimen1 by \dimen7
\advance\dimen2 by \dimen8
\multiply\dimen4 by -#3
\rlap{ \kern \dimen1 \raise \dimen2 
\hbox{ \li {\dimen4} {\dimen5} 0pt }}
}
\fi}


\def\wr #1 #2 #3 {
{
\multiply\dimen1 by #1
\advance\dimen1 by \dimen7
\multiply\dimen2 by -#2
\advance\dimen2 by \dimen8
\rlap{ \kern \dimen1 \raise \dimen2 \hbox{\tenrm #3 }}}}


\def\wrInBox #1 #2 #3
{
\dimen9=\dimen1
\multiply\dimen9 by #1
\advance\dimen9 by \dimen7
\dimen10=\dimen1
\multiply\dimen10 by -#2
\advance\dimen10 by \dimen8
\rlap{ \kern \dimen9  \raise \dimen10
\vbox to \dimen2{
\vfill \hbox to \dimen1{\hfill {\tenrm #3} \hfill } \vfill }}
}

\def\SpecPicture#1#2#3\endpicture
   {\multiply\dimen0 by #1
    \multiply\dimen6 by #2
    \vbox to 0pt 
   {\vss \vskip \dimen6 \hbox to \dimen0{\nullfont #3 \hss}\vss}
   }

\def\Picture#1#2\endpicture
   {\multiply\dimen0 by #1 
    \bigskip \centerline{\hbox to \dimen0 
   {\nullfont #2 \hss}} \nobreak \bigskip}

\def\CapPicture#1#2#3\endpicture
   {\multiply\dimen0 by #1 
    \bigskip \centerline{\hbox to \dimen0 
   {\nullfont #3 \hss}} \nobreak \bigskip \centerline{#2} \bigskip}

\NoBlackBoxes

\def\l{\lambda}
\def\Pm{P^*_\mu}
\def\Pl{P_\l}
\def\Qm{Q_\mu}
\def\t{\theta}
\def\a{\alpha}
\def\b{\beta}
\def\g{\gamma}
\def\d{\delta}
\def\L{\Lambda}
\def\Lt{\L^\t}
\def\C{\Bbb C}
\def\R{\Bbb R}
\def\Z{\Bbb Z}
\def\Q{\Bbb Q}
\def\T{\Bbb T}

\def\tht{\thetag}

\def\plm{\psi_{\l/\mu}}
\def\G{\Gamma}
\def\ap{\a^+}
\def\am{\a^-}
\def\api{\ap_i}
\def\ami{\am_i}
\def\bp{\b^+}
\def\bm{\b^-}
\def\bpi{\bp_i}
\def\bmi{\bm_i}
\def\lp{\l^+}
\def\lm{\l^-}

\def\gp{\g^+}
\def\gm{\g^-}
\def\ABG{\left(\matrix \ap;\bp;\gp \\ 
\am;\bm;\gm \endmatrix\right)}
\def\pin{\phi_\infty}
\def\psin{\psi_\infty}
\def\tri{\therosteritem}
\def\e#1{e^{2\pi i#1}}
\def\r{\rho}
\def\N#1{{\frak N}(#1)}
\def\ep{\varepsilon}
\def\dd#1{\dot{#1}}
\def\ddd#1{\ddot{#1}}
\def\Ti{\T^\infty}
\def\P{\Phi}
\def\U{\Upsilon^\t}
\def\Un{\Upsilon^\t_n}
\def\Unm{\Upsilon^\t_{n-1}}
\def\Ex{\operatorname{Ex}}
\def\tr{\operatorname{tr}}

\def\Zgln{{\frak Z}(\frak{gl}(n))}

\def\Sm{{\Bbb S}_\mu}
\def\Zg{{\frak Z}(\frak{g})}
\def\Ug{{\Cal U}(\frak{g})}
\def\ite#1{\item"(#1)"}
\def\Id{\operatorname{Id}}

\def\A{\frak{A}}
\def\Ab{\overline{\A}}
\def\B{\frak{B}}
\def\fC{\frak{C}}
\def\om{\omega}
\def\Proj#1#2{\operatorname{Proj}^{#1}_{#2}}
\def\bProj#1#2{\overline{\operatorname{Proj}}^{#1}_{#2}}

\rightheadtext{Asymptotics of Jack polynomials}
\leftheadtext{A.~Okounkov and G.~Olshanski}

\topmatter
\title
Asymptotics of Jack polynomials as the number of
variables goes to infinity 
\endtitle
\author
Andrei Okounkov and Grigori Olshanski
\endauthor 
\abstract
In this paper we study the asymptotic behavior of the Jack rational 
functions $P_\l(z_1,\dots,z_n;\t)$ as the number of variables 
$n$ and the signature $\l$ grow to infinity. Our results 
generalize the results of A.~Vershik and S.~Kerov \cite{VK2}
obtained in the Schur function case ($\t=1$). For $\t=1/2,2$
our results describe approximation of the spherical functions
of the infinite-dimensional symmetric spaces $U(\infty)/O(\infty)$
and $U(2\infty)/Sp(\infty)$ by the spherical functions of
the corresponding finite-dimensional symmetric spaces. 
\endabstract
\address 
A.O.:  Dept.\ of Mathematics, University of
Chicago, 5734 South University Ave., Chicago, IL 60637-1546.
{\it E-mail address:\/} {\rm okounkov\@math.uchicago.edu}
\newline\vskip 0.1 in 
G.O.: Institute for Problems of Information Transmission,
Bolshoy Karetny 19, 101447 Moscow GSP--4, Russia.
{\it E-mail address:\/} {\rm olsh\@ippi.ac.msk.su}
\endaddress
\thanks 
The authors were supported by the Russian Basic Research Foundation
grant 95-01-00814. The first author's stay at IAS in Princeton and
MSRI in Berkeley was supported by  NSF grants 
DMS--9304580 and DMS--9022140 respectively.
\endthanks 
\toc
\head 1. Introduction\endhead
\subhead 1.1. Statement of the main result \endsubhead
\subhead 1.2. Regular and infinitesimally regular sequences\endsubhead
\subhead 1.3. Extremality of the limit functions \endsubhead
\subhead 1.4. Related results \endsubhead
\subhead 1.5. Acknowledgments \endsubhead
\head 2. Jack polynomials and shifted Jack polynomials
\endhead
\subhead 2.1. Orthogonality
\endsubhead 
\subhead 2.2. Interpolation
\endsubhead 
\subhead 2.3. Branching rules
\endsubhead
\subhead 2.4. Binomial formula 
\endsubhead 
\subhead 2.5. Generating functions 
\endsubhead 
\subhead 2.6. Partitions and signatures  
\endsubhead 
\subhead 2.7. Extended symmetric functions  
\endsubhead 
\head
3. Asymptotic properties of Vershik-Kerov sequences of signatures
\endhead
\head
4. Sufficient conditions of regularity
\endhead
\head
5. Necessary conditions of regularity
\endhead
\subhead 
5.1. The ``only if'' part of Theorem 1.1
\endsubhead
\subhead 
5.2. A growth estimate for $|f(\l)|$, $f\in\Lt$
\endsubhead
\head
6. Proof of Theorem 1.3
\endhead
\head 7. Appendix. A direct proof of the formula \tht{2.10}
for generating functions
\endhead
\endtoc
\endtopmatter

\newpage 

\document

\head
1.~Introduction
\endhead

\subhead
1.1 Statement of the main result
\endsubhead 
Jack symmetric functions
$$
\Pl(z_1,\dots,z_n;\t)\in\Q(\t)[z^{\pm1}]^{S(n)}
$$
which are indexed by decreasing sequences of integers
(called signatures)
$$
\l=(\l_1\ge\dots\ge\l_n)\in\Z^n\,,
$$
are eigenfunctions of the quantum Calogero-Sutherland
Hamiltonian \cite{C,Su}
$$
H_{CS}=\frac12\sum_i 
\left(z_i \frac{\partial}{\partial z_i}\right)^2 +
\t \sum_{i\ne j} \frac{z_i^2}{z_i-z_j} 
\frac{\partial}{\partial z_i}\,. \tag 1.1
$$
The Hamiltonian $H_{CS}$ describes $n$ identical quantum 
particles living on the circle
$\T=\{z\in\C, |z|=1\}$ and interacting via an inverse-square
potential. Here $\t>0$ is a fixed positive parameter
(a coupling constant). For 3
special values
$$
\t=1,\frac12,2
$$
the functions $\Pl(z;\t)$ reduce to Schur functions
(which are spherical functions of the symmetric space
$U(n)$) and spherical functions of Gelfand pairs
$$ 
U(n)\supset O(n)\,, \quad U(2n)\supset Sp(n) \tag 1.2 
$$
respectively. For these values of $\t$,  the meaning of the 
parameter $\t$ is half the multiplicity of a root
in the restricted root
system (of type $A_n$) associated to the corresponding
symmetric space.

In this paper we study the behavior of the Jack polynomials
as the number of variables grows to infinity in such a way that
all but finitely many variables $z_i$ are kept equal to 1. 
More precisely, set $\T^k=(\T)^{\times k}$ and consider the
direct limit 
$$
\Ti:=\varinjlim\T^k
$$ 
with respect to embeddings $(z_1,\dots,z_k)\mapsto(z_1,\dots,z_k,1)$.
Equip $\Ti$ with the direct limit topology. 

Let $\l(n)=(\l_1(n)\ge\dots\ge\l_n(n))$ be a sequence of signatures.
We call this sequence {\it regular} if the following sequence of
{\it normalized Jack functions}
$$
\P_{\l(n)}(z_1,\dots,z_n;\t):=
\frac{P_{\l(n)}(z_1,\dots,z_n;\t)}
{P_{\l(n)}(1,\dots,1;\t)} \tag 1.3 
$$
converges uniformly on the compact subsets of $\Ti$ (it suffices to
consider the subsets $\T^k\subset\Ti$). Notice that the normalization
in \tht{1.3} is necessary to make these functions converge at least at
the point $(1,1,\dots)\in\Ti$. Note also that the normalization \tht{1.3}
is precisely the normalization for the spherical function (for 
$\t=1,\frac12,2$). 

Our main result is an explicit description of all regular sequences
and an explicit computation of the corresponding limits 
$$
\lim_{n\to\infty} \P_{\l(n)}(z;\t)\,, \tag 1.4 
$$
given in  Theorem 1.1 below. In the Schur function case ($\t=1$) this
program was carried out by A.~Vershik and S.~Kerov
\footnote{
we order the names of A.~Vershik and S.~Kerov in the Russian
alphabetical order}
(see \cite{VK2})
in their proof of the classification theorem for characters of
the group $U(\infty)$.   
In particular they have found (for $\t=1$) 
a beautiful description of the set of
regular  sequences. Their necessary and sufficient
conditions, which we shall call the{\it Vershik-Kerov conditions},
are as follows. 

First, suppose we have a sequence of partitions 
$
\left\{\l(n)\right\}_{n=1,2,\dots}
$\,.
This sequence is said to be a {\it Vershik-Kerov sequence}
(or VK for short) if
the following limits exist:
$$
\alignat2
\a_i &:= \lim \frac{\l(n)_i}n &\quad&<\infty\,,\\
\b_i &:= \lim \frac{\l'(n)_i}n &\quad&<\infty\,,\\
\d &:= \lim \frac{|\l(n)|}n &\quad&<\infty\,.\\
\endalignat
$$
Here $\l'$ denotes the partition conjugated to $\l$; in
other words, $\l'_i$ is the length of the $i$-th column in the
diagram of $\l$. It is easy to see that the number
$$
\g:=\d - \sum (\a_i + \b_i) \ge 0
$$
is nonnegative. The numbers $\a_i$, $\b_i$, $\g$ are called the
{\it VK parameters} of the sequence $\l(n)$.

Now, with each signature $\l$ one can associate two partitions $\lp$
and $\lm$ as follows. Suppose $\l_p\ge 0\ge \l_{p+1}$ for some
$p=1,\dots,n$. Then, by definition,
$$
\align
\lp&=(\l_1,\dots,\l_p)\,,\\
\lm&=(-\l_n,\dots,-\l_{p+1})\,.
\endalign
$$
A sequence of signatures 
$$
\l(n)=(\l_1(n)\ge \dots \ge \l_n(n))
$$ 
is called a {\it VK sequence of signatures} if the 
two sequences of partitions $\lp(n)$ and $\lm(n)$ are
Vershik-Kerov sequences. The corresponding numbers
$\a^\pm_i$, $\b^\pm_i$, and $\g^\pm$ are called the
VK parameters of the sequence $\l(n)$. A typical member of
a typical VK sequence looks like this: 
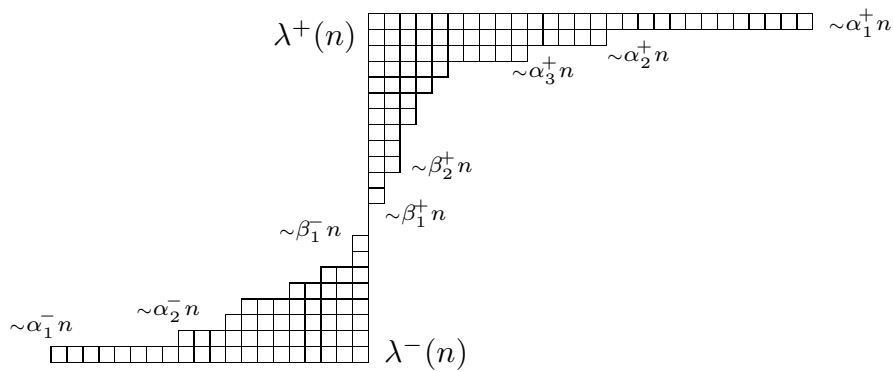
\CapPicture{58}{Fig.~1 A typical element of a Vershik-Kerov sequence}{
\SetGrid{5 pt}{5 pt}
\LineWidth{0.1 pt}
\SetOrigin -2 0  
\hl 0 0 28
\hl 0 1 28
\hl 0 2 15
\hl 0 3 10
\hl 0 4 5
\hl 0 5 4
\hl 0 6 3
\hl 0 7 3
\hl 0 8 2
\hl 0 9 2
\hl 0 10 2
\hl 0 11 1
\hl 0 12 1
\vl 0 0 22
\vl 1 0 12
\vl 2 0 10
\vl 3 0 7
\vl 4 0 5
\vl 5 0 4
\vl 6 0 3
\vl 7 0 3
\vl 8 0 3
\vl 9 0 3
\vl 10 0 3
\vl 11 0 2
\vl 12 0 2
\vl 13 0 2
\vl 14 0 2
\vl 15 0 2
\vl 16 0 1
\vl 17 0 1
\vl 18 0 1
\vl 19 0 1
\vl 20 0 1
\vl 21 0 1
\vl 22 0 1
\vl 23 0 1
\vl 24 0 1
\vl 25 0 1
\vl 26 0 1
\vl 27 0 1
\vl 28 0 1
\wr -6 2 {$\lambda^+(n)$}
\wr 29 1 {$\sssize{\sim}\ssize{\a^{\!+}_1 n}$}
\wr 15 3 {$\sssize{\sim}\ssize{\a^{\!+}_2 n}$}
\wr 9 4 {$\sssize{\sim}\ssize{\a^{\!+}_3 n}$}
\wr 1 13 {$\sssize{\sim}\ssize{\b^{\!+}_1 n}$}
\wr 2 10.5 {\,\,$\sssize{\sim}\ssize{\b^{\!+}_2 n}$}
\SetOrigin -2  22 
\hl 0 0 -20
\hl 0 -1 -20
\hl 0 -2 -12
\hl 0 -3 -9
\hl 0 -4 -8
\hl 0 -5 -5
\hl 0 -6 -3
\hl 0 -7 -1
\hl 0 -8 -1
\vl -1 0 -8
\vl -2 0 -6
\vl -3 0 -6
\vl -4 0 -5
\vl -5 0 -5
\vl -6 0 -4
\vl -7 0 -4
\vl -8 0 -4
\vl -9 0 -3
\vl -10 0 -2
\vl -11 0 -2
\vl -12 0 -2
\vl -13 0 -1
\vl -14 0 -1
\vl -15 0 -1
\vl -16 0 -1
\vl -17 0 -1
\vl -18 0 -1
\vl -19 0 -1
\vl -20 0 -1
\wr 1 0 {$\lambda^-(n)$}
\wr -23 -2 {$\,\sssize{\sim}\ssize{\a^{\!-}_1 n}$}
\wr -15 -3 {$\,\sssize{\sim}\ssize{\a^{\!-}_2 n}$}
\wr -6 -8 {$\,\sssize{\sim}\ssize{\b^{\!-}_1 n}$}
}
\endpicture
Throughout this paper we assume the following important 
notational convention. For a sequence $(\l(n))$ of signatures the 
letter $n$ will always stand for the length of the 
signature $\l(n)$. The values of $n$ may range, however,
over some proper infinite subset of $\{1,2,\dots\}$.
In other words, we allow gaps in the values of $n$. 

Our main result, which is a direct generalization of the result
of Vershik and Kerov \cite{VK2} is the following: 

\proclaim{Theorem 1.1} For any $\t>0$, a sequence $\l(n)$ of signatures is
regular if and only if it is VK. If $\l(n)$ is VK with
parameters $\a^\pm_i$, $\b^\pm_i$, and $\g^\pm$ then 
we have 
$$
\lim_{n\to\infty} \P_{\l(n)}(z;\t)= \phi_{\a,\b,\g}(z_1) 
\phi_{\a,\b,\g}(z_2) \dots\,, \tag 1.5
$$
where the function $\phi_{\a,\b,\g}(z)$ is the following product 
$$
\phi_{\a,\b,\g}(z):=e^{\gp (z-1) + \gm (z^{-1}-1)} \prod_i 
\frac{(1+ \bpi (z-1))}
{(1- \api (z-1)/ \t)^{\t}} 
\frac{(1+\bmi (z^{-1}-1))}
{(1-\ami  (z^{-1}-1)/ \t)^{\t}}\,. \tag 1.6
$$
\endproclaim

The ``if'' part of this theorem will be established in Theorem 4.1 ;
the ``only if'' part will follow from Theorem 5.1 in Section 5.1.

We shall denote by $\P_{\a,\b,\g}(z)$ the function in the right-hand
side of \tht{1.5}.
The  remarkable twofold factorization of the limit $\P_{\a,\b,\g}$ has
the following interpretation  in the group-theoretic situation 
$\t=1,\frac12,2$ (see \cite{Vo,Ol1}). The factorization \tht{1.5}, which states
that all variables $z_i$ become, in a sense, independent of each
other, is the infinite-dimensional
degeneration of the functional equation satisfied by the
spherical functions. It reflects the fact that the 
convolution algebra of $K(n)$-biinvariant functions on $G(n)$ 
(here $G(n)\supset K(n)$ stands for a series of Gelfand pairs)
becomes, in a certain precise sense, a semigroup algebra as $n\to\infty$.  

The second factorization \tht{1.6} asserts that
all spherical representations of the corresponding Gelfand pairs
are tensor products of elementary (anti)-bosonic, (anti)-fermionic,
and certain intermediate representations. 

In the context of the quantum
Calogero-Sutherland model the parameters $\a$ and $\b$ are 
(in a somewhat different limit transition) associated with
quasi-particles and 
quasi-holes, see \cite{Ha,LPS,Ok11}.  

\subhead
1.2 Regular and infinitesimally regular sequences
\endsubhead

Although there are several more or less explicit formulas for
Jack polynomials (see, for example, \tht{2.3} below),
 none of them seems to be suitable for direct
investigation of the limit \tht{1.4}. The strategy we shall use
in this paper (and which in an implicit form was present already
in the work of Vershik and Kerov) is to replace the study of 
the uniform convergence of \tht{1.3} by the investigation of 
convergence of the coefficients in the Taylor expansion 
of \tht{1.3} about the point $z=(1,1,\dots)$. 
The uniform convergence of \tht{1.3} is, in fact, equivalent
to weak convergence of certain measures (see below); the study of 
its Taylor series is then equivalent to the study of
moments of those measures, which is a time-honored way of
proving limit theorems. In the group-theoretic
situation $\t=1,\frac12,2$, this means that one replaces the
study of functions on a Lie group $G$ by the study of the
corresponding linear functionals  on the universal enveloping algebra
of the Lie algebra $\frak g$ of $G$. 

Let us call a sequence $\l(n)$ an {\it infinitesimally
regular} sequence if the Taylor series of the functions \tht{1.3}
about the point $(1,1,\dots)$ tend to a limit  coefficient by 
coefficient. Since it is very often that quite complicated
special functions have rather simple Taylor series, it is no
surprise that this condition turns out to be  easier to verify
than the regularity condition. Of course, this condition
is not {\it a priori} equivalent to the regularity 
condition. However we shall prove the following result
(the last condition \tri{iv} in that theorem will be
explained below)

\proclaim{Theorem 1.2} Given a sequence $\l(n)$, 
the following conditions are equivalent
\roster
\item"(i)" $\l(n)$ is regular,
\item"(ii)" $\l(n)$ is infinitesimally regular,
\item"(iii)" $\l(n)$ is a Vershik-Kerov sequence,
\item"(iv)" for any polynomial $f\in\Lt$ the limit $\lim_{n\to\infty} 
\frac{\tsize f(\l(n))}{\tsize n^{\deg f}}$ exists,
\endroster
and the 4 above conditions are also equivalent to the hypothesis of
Theorem 5.1.
\endproclaim

The proof of this theorem will be completed in Section 5.1. 

The meaning of the condition \tri{iv} is the following. The 
quantum Calogero-Sutherland Hamiltonian \tht{1.1} is completely
integrable: there are $n$ independent commuting differential
operators which commute with $H_{CS}$, or, in other words,
there are $n$ independent quantum integrals of motion, see e.g.\
\cite{He}. The values
of those conserving quantities in the quantum state described by $P_\l(z)$
are polynomials in $\l_1,\dots,\l_n$ symmetric in variables
$\l_i-\t\,i$. Conversely, each polynomial with such symmetry
corresponds to some integral of motion. Thus, the subalgebra
$$
\Lt(n)\subset\C[\l_1,\dots,\l_n]
$$ 
of polynomials symmetric in variables $\l_i-\t i$ is a
natural commutative algebra of observables in the model. For example,
the polynomial
$$
\l_1+\dots+\l_n \in \Lt(n)
$$
corresponds to the differential operator 
$\sum z_i\frac{\partial}{\partial z_i}$ which is the total momentum
of the system. The polynomial
$$
E(\l)=\frac 12\sum_i \l_i^2 + \t \sum_i (n-i)\l_i
$$
corresponds to the energy of the system, that is, to the operator
$H_{CS}$ itself.

We denote by
$$
\Lt:=\varprojlim \Lt(n),
$$
the inverse limit of $\Lt(n)$ as filtered (by degree of polynomials)
algebras with respect to homomorphisms:
$$
\align
\Lt(n+1)\quad  &\to  \quad \Lt(n) \,, \\
f(\l_1,\dots,\l_{n+1})&\mapsto f(\l_1,\dots,\l_n,0)\,. \tag 1.7
\endalign
$$
Examples of elements of $\Lt$ are the functions $\sum\l_i$,
$\sum \l_i^2 - 2\t\sum i \l_i$, or, more generally,
$\sum ((\l_i - \t i)^k-(-\t i)^k)$, where $k=0,1,2,\dots$.    

The condition \tri{iv} now means that all ``stable'' observables $f$,
after proper normalization, have a limit value as  $n\to\infty$. 
Observe that energy 
$$
E(\l)=
\frac 12\left(\sum_i \l_i^2 -  2\t \sum_i i\l_i \right) + 
n\t\sum_i \l_i 
$$
is not a stable polynomial. Notice, however, that 
if the condition \tri{iv} is satisfied then the limit
$$
\lim_{n\to\infty} \frac{E(\l)}{n^2}
$$
exists. More generally, one can consider the algebra $\Lt[n]$
of polynomials in $n$ with coefficients in stable polynomials
in $\l$ and introduce the notion of degree in this algebra by setting
$$
\deg n = \deg \l_i = 1 \,.
$$
Then, for example, $E(\l)$ is an element of this algebra of degree 2.
It is clear that \tri{iv} is equivalent to the more general 
condition
$$
\forall f\in\Lt[n]\quad \exists \lim_{n\to\infty} 
\frac{f(\l,n)}{n^{\deg f}} \,.
$$

Note finally that, although we shall make no use of those $n$ commuting
differential operators in the present paper, their role is
crucial in the proof of the binomial theorem (see Section 2.4). 

\subhead
1.3 Extremality of the limit functions
\endsubhead

For $n=1,2,\dots$ define the following convex subsets 
$\Un\subset C(\T^n)$. The subset $\Un$ consists, by definition, of the
functions of the form
$$
\phi(z_1,\dots,z_n)=\sum_\l c_\l \P_\l(z;\t)\,,
$$
where $\l$ ranges over signatures of length $n$ and 
the coefficients $c_\l$ are arbitrary satisfying 
$$
c_\l \ge 0\,, \quad \sum_\l c_\l = 1\,.
$$
It is clear that $\Un$ is a simplex with countably many
vertices. It follows from the positivity of the branching coefficients 
 for Jack polynomials (see the formula \tht{2.3} below) that
the simplex $\Un$ is contained in (and for $n=1$ coincides with)
the simplex of characteristic
functions of probability measures on $\Z^n$. 

It also follows from the same positivity  that the map
$$
\phi(z_1,\dots,z_n)\mapsto\phi(z_1,\dots,z_{n-1},1)
$$
defines an affine map 
$$
\Un \to \Unm \,.
$$
We set 
$$
\U:=\varprojlim \Un\,.
$$
By definition, each element of $\U$ is a function on $\Ti$.
It is clear that the set $\U$ is convex. Denote by 
$\Ex \U$ the set of its extreme points. Recall that
$\P_{\a,\b,\g}(z)$, $z\in\Ti$, denotes the RHS of the formula
\tht{1.5}. In Section 6 we
shall demonstrate the following

\proclaim{Theorem 1.3} We have
$$
\Ex \U = \big\{\P_{\a,\b,\g}\big\}_{\a^\pm,\b^\pm,\g^\pm}\,, 
$$
where $\a^\pm,\b^\pm,\g^\pm$ range over all possible values of VK 
parameters:
$$
\gather
\ap_1 \ge \ap_2 \ge \dots 0\,, \quad
\am_1 \ge \am_2 \ge \dots 0\,,  
\\
\bp_1 \ge \bp_2 \ge \dots 0\,, \quad
\bm_1 \ge \bm_2 \ge \dots 0\,, \\
\sum (\api+\ami + \bpi+\bmi) < \infty\,,\\
\b^+_1 + \b^-_1 \le 1\,,\\
\g^\pm \ge 0\,.
\endgather
$$
\endproclaim

Theorem 1.3 is a rather
straightforward corollary of Theorem 1.1 and some general abstract
theorems, see Section 6.

The statement of Theorem 1.3 was known for $\t=1,\frac12,2$ where
the set $\Ex \U$ is the set of characters of 
$$
U(\infty)=\varinjlim U(n)
$$ 
and the set of spherical functions of Gelfand pairs 
$$
U(\infty)\supset O(\infty)\,, \quad U(2\infty)\supset Sp(\infty) \tag 1.8
$$
respectively. 

The description of $\Upsilon^1$ is a famous theorem,
various pieces of which were proved by A.~Edrei, D.~Voiculescu,
R.~Boyer, and  A.~Vershik and S.~Kerov. Namely, D.~Voiculescu 
constructed in \cite{Vo} a large supply of characters of $U(\infty)$
and conjectured that his list is complete. Later R.~Boyer \cite{B1}
and independently A.~Vershik and S.~Kerov \cite{VK2} observed that that
completeness follows from an old hard theorem of A.~Edrei
\cite{Ed}. In the same note, A.~Vershik and S.~Kerov have also outlined an
alternative  plan of proof (which uses approximation of
characters of $U(\infty)$ by characters of $U(n)$) and have 
formulated the $\t=1$ case of above Theorem 1.1. 
A complete proof of a particular case of that result of
Vershik and Kerov (the case when $\l(n)$ is assumed to be a partition)
was published later by R.~Boyer in \cite{B2}. In that proof he used
techniques different from those announced by Vershik and Kerov (and 
different from those used in our paper). 

The description of $\U$ for $\t=\frac12,2$ was obtained 
in the papers  of  G.~Olshanski
and D.~Pickrell. Spherical (and more general admissible)
representations of the infinite-dimensional Gelfand pairs
corresponding to classical series of Riemannian symmetric spaces
were constructed and studied by Olshanski in \cite{Ol1}.
The proof of the completeness of Olshanski's lists of spherical
representations was completed  by D.~Pickrell
in \cite{P} by using embeddings of symmetric spaces into each other. 

The new piece of information we obtain in the case
$\t=\frac12,2$ is how exactly the spherical functions of the
Gelfand pairs \tht{1.2} approximate the spherical functions of the
corresponding infinite-dimensional Gelfand pairs \tht{1.8}. 

\subhead
1.4~Related results
\endsubhead
\subhead
1.4.1~Shifted Schur functions and binomial formula 
\endsubhead

Our research presented in this paper started in 1994 in the
framework of the general classification problem of 
admissible representations of infinite-dimensional Gelfand 
pairs, see \cite{Ol1-2,KO,OV,Ok1-2}. We began by analyzing the
details of the proofs sketched in the seminal paper \cite{VK2} of
A.~Vershik and S.~Kerov. As we already mentioned above, that
paper contained, among other things, the following important idea:
one has to consider the convergence of Taylor series coefficients.

The Taylor expansion of the Schur functions about the point
$(1,\dots,1)$
are given by the following {\it binomial formula}
$$
\frac{s_\l(1+x_1,\dots,1+x_n)}{s_\l(1,\dots,1)}
=\sum_\mu \frac{s^*_\mu(\l) \, s_\mu(x)}
{\prod_{(i,j)\in\mu}(n+j-i)}\,. \tag 1.9
$$
Here $s^*_\mu$ are certain factorial analogs of Schur
functions which can be defined by the following determinant
ratio formula
$$
s^*_\mu(x_1,\dots,x_n)=
\frac{\det
\big[(x_i+n-i)\cdots(x_i-i+j-\mu_j+1)\big]_{1\le i,j\le n}}
{\prod_{1\le i<j\le n} (x_i-x_j+j-i)} \,. \tag 1.10
$$
We call these functions the {\it shifted Schur functions}, see
\cite{OO}. They
differ by a shift of variables only from the factorial Schur
functions introduced in \cite{BL} and studied further in \cite{M2}
and other papers cited in \cite{OO}. 

A formula equivalent to the  binomial formula \tht{1.9}
was established by A.~Lascoux, see \cite{Lasc} and also Example I.3.10
in \cite{M2}. Its proof is elementary and
based on the determinant ratio formula
for the Schur functions
$$
s_\mu(x_1,\dots,x_n)=
\frac{\det
\big[x_i^{\mu_j+n-j}\big]_{1\le i,j\le n}}
{\prod_{1\le i<j\le n} (x_i-x_j)} \,. \tag 1.11
$$
The observation that the coefficients appearing in Lascoux's binomial 
formula  are factorial analogs of ordinary Schur functions
was made by Olshanski in \cite{Ol3}, see also \cite{OO} and the next
subsection. 

{F}rom \tht{1.9} it is clear that for $\t=1$ a sequence $\l(n)$ of signatures is 
infinitesimally regular if and only if the limit
$$
\lim_{n\to\infty} \frac{s^*_\mu(\l(n))}{n^{|\mu|}}
$$
exists for every partition $\mu$. 

\subhead
1.4.2~Quantum immanants(\cite{OO,Ok6-7})
\endsubhead

For general Jack polynomials (or even for $\t=1/2,2$), no
formulas as simple as \tht{1.11} are available. Therefore, in order to
generalize \tht{1.9}, one needs some insight into the invariant 
meaning of the expansion \tht{1.9}.

Observe that in the binomial formula \tht{1.9} there is a certain 
symmetry between $\l$ and $x$. There exists, in fact, a very precise
parallel between the functions $s_\mu$ and $s^*_\mu$. 
In short, the relation between them is the same
as the relation between the algebra 
$$
Z(G)=C(G)^G \tag 1.12
$$
of central continuous functions on a compact group $G$ and 
the algebra of the {\it Laplace operators}
$$
\Zg=\Ug^G \tag 1.13 
$$
on the group $G$. Here $\frak{g}$ is the (complexified) Lie 
algebra of $G$ and $\Ug$ stands for its universal enveloping algebra.
(In our case, $G=U(n)$). Instead of the whole group  $G$,  one
can look at its maximal torus and the radial parts of the elements of $Z(G)$
and $\Zg$. 

The functions $s_\mu$ and $s^*_\mu$ correspond to
distinguished linear bases of the spaces \tht{1.12} and \tht{1.13}
respectively. In \tht{1.12}, the Haar measure defines a natural
scalar product
$$
(f_1,f_2)=\int_G f_1(g) \, \overline{f_2(g)} \, dg\,,
$$
with respect to which the characters (i.e.\ Schur functions)
form an orthogonal basis. Another way of expressing this orthogonality
is to say that a character (as a function on $G$) vanishes in all
but one irreducible representations. 

In the algebra $\Zg$, there is no natural scalar product. However,
there is a natural pairing between \tht{1.12} and \tht{1.13}, namely
$$
\langle D,f\rangle = \big(Df\big) (1,\dots,1)\,, \quad D\in\Zg\,,
$$
where $D$ is the radial part of a Laplace operator and $f$ is
the radial part of an element of $Z(G)$. In particular,
$$
\langle D,s_\l\rangle =\tr \pi_\l(D)\,,
$$
where $\pi_\l$ stands for the irreducible representation with
highest weight $\l$. There exists a linear basis  
$\left\{\Sm
\right\} \subset \Zg$ indexed by partitions $\mu$ with at most 
$n$ parts such that
\roster
\item $\Sm$ is an operator of order $|\mu|$\,,
\item $\pi_\l(\Sm)=0$ for all partitions $\l$ such that
$\mu\not\subset\l$\,.
\endroster 
In words, $\Sm$ vanish in as many 
irreducible representations as possible (it is clearly impossible
for an element of $\Zg$ to vanish in all but one irreducible 
representation). 
We call these elements of $\Zg$ the {\it quantum immanants},
see \cite{OO,Ok6-7}. Their
relation to the $s^*$-functions is the following:
$$
\pi_\l({\Bbb S}_\mu)=s^*_\mu(\l) \cdot \Id\,,
$$
which follows easily from the definition of $\Sm$, the
Harish-Chandra isomorphism theorem, and the following three 
direct consequences of \tht{1.10}:
$$
\align
&\deg s^*_\mu=|\mu|\,,
\tag 1.14 \\
&\text{
the polynomial $s^*_\mu(\l)$ is symmetric in variables $\l_i-i$}\,,
\tag 1.15 \\
&\text{
$s^*_\mu(\l)=0$ for all partitions $\l$ such that
$\mu\not\subset\l$}\,.
\tag 1.16
\endalign
$$

The quantum immanants and shifted Schur functions enjoy  quite
a few remarkable properties,  some of which do and some of which 
don't have analogs for the ordinary Schur functions. We studied
them at length in our papers \cite{OO,Ok6-7}. For example, the
explicit formulas for quantum immanants proved in \cite{Ok6}
(see also \cite{N2,Ok7}) generalize the classical Capelli identities
(see, for example, the papers \cite{HU,N1} and also \cite{S1}). 

It is amazing how powerful is the vanishing property \tht{1.16}.
For example, it yields a short 
computation-less proof of the expansion \tht{1.9}
(see \cite{OO}, Section 5) and also other previously known properties
of the $s^*$-functions.

To summarize, the role which the quantum immanants 
play behind the scene in the Vershik-Kerov theorem is the 
following. When one goes from the Schur functions to their Taylor
expansions, one replaces the study of the functions 
$$
\alignat3
&\frac{\tr\pi_\l(g)}{\dim \pi_\l}\,, \quad &&g\in U(n)\,,&& 
\quad \pi_\l\in U(N)^\wedge\,, \quad n\le N \,,  \\
\intertext{by the study of the functions}
&\frac{\tr \pi_\l(X)}{\dim \pi_\l}\,, \quad &&X\in\Zgln\,,&&
\quad \pi_\l\in U(N)^\wedge\,, \quad n\le N \,. \tag 1.17
\endalignat
$$
The quantum immanants are the distinguished elements of $\Zgln$
on that the functions \tht{1.17} are easy to compute (see the formulas
in \cite{OO}, Section 10, and also in \cite{Ok6}, Section 5). 

\subhead
1.4.3~The case of general $\t>0$ 
\endsubhead

For general $\t>0$, we have the algebra of quantum integrals of motion in
the Calogero-Sutherland model instead of the algebra of the radial parts of Laplace
operators. Instead of quantum immanants, we have differential 
operators that annihilate as many Jack polynomials $P_\l(x;\t)$
as possible. {F}rom the discussion after Theorem 1.2 it is clear that
the role of the $s^*_\mu$ functions is to be played by the polynomials
$\Pm(x;\t)$ satisfying the following conditions (compare with \tht{1.14-16}):
$$
\align
&\deg \Pm=|\mu|\,,
\tag 1.18 \\
&\text{
the polynomial $\Pm(\l;\t)$ is symmetric in variables $\l_i-\t i$}\,,
\tag 1.19 \\
&\text{
$\Pm(\l;\t)=0$ for all partitions $\l$ such that
$\mu\not\subset\l$}\,.
\tag 1.20 
\endalign
$$
Note that the polynomials $\Pm(x;\t)$ are a particular case of the polynomials 
considered by S.~Sahi in \cite{S1}. We call these polynomials 
the {\it shifted Jack polynomials}; they shall play a central role in 
the present paper and will be discussed more formally in Section 2. 

Modifying suitably (see \cite{OO2})
the argument from  \cite{OO}, Section 5 (2nd proof),
one deduces from the definition
 \tht{1.18-20} the binomial formula given in the Section 2.4 below.
However, the abstract formula \tht{2.6} alone is not sufficient for the purposes 
of the present paper. One needs some control of the polynomials 
$\Pm(x;\t)$, for example, one would like  to know their highest degree
term. We have conjectured the explicit formula \tht{2.4} for 
the shifted Jack polynomials $\Pm(x;\t)$, which, in particular, implies
that 
$$
\Pm(x;\t)=P_\mu(x;\t)+\dots\,, \tag 1.21
$$
where dots stand for lower degree terms. However, we were only able
to prove our conjecture  for $\mu=(m)$ (using the argument reproduced
in Section 7) and to check \tht{1.21} in the group-theoretic situation  $\t=\frac12,2$. 

First proof of \tht{1.21} for general $\t$ 
was given, together with some other fundamental results, 
by F.~Knop and S.~Sahi in \cite{KS}; see also
their papers \cite{Kn,S}. They used certain difference equations
for these polynomials. 
\footnote{It seems to be an interesting and important problem to
understand better the nature of the Knop-Sahi difference equations.
For example, their analogs do not exists for other series of root
systems, see \cite{Ok10}.
}
\footnote
{
Also note, that by virtue of the binomial formula \tht{2.6}, the 
relation \tht{1.21} is equivalent to formula \tht{4.2} in \cite{OO2}
for the Bessel functions. Apparently, that formula for the
Bessel functions was also known
to M.~Lassalle, see the paper \cite{BF}, Sections 3 and 6, where the authors 
cite Lassalle's unpublished letters. Lassalle's argument is very different
from proofs of \tht{1.21} given in \cite{KS} and \cite{Ok8}.   
} 

Later the formula \tht{2.4} and analogs of
other properties of shifted Schur functions were established for
shifted Jack polynomials and also their $q$-analogs in \cite{OO2,Ok8-9}.

\subhead
1.4.4~Log-concavity and the ``only if'' part of Theorem 1.1 
\endsubhead

The most subtle point of the ``only if'' part of Theorem 1.1
seems to be to justify the following implication 
(see the statement of Theorem 1.2)
$$
\text{\tri{i}}\overset?\to\Longrightarrow
\text{\tri{ii}} \,. \tag 1.22
$$
This point is quite subtle even in the $\t=1$ case.
(A.~Vershik and S.~Kerov did not discuss the proof of the
``only if'' part of their theorem in \cite{VK2}.)

A proof and an explanation  of \tht{1.22} based on some general log-concavity
results was proposed in the paper \cite{Ok3}. (See also the subsequent 
papers \cite{Gr,Ok4-5,Ka}.) Similar log-concavity argument applies in
the situation considered in \cite{OV}.  

Unfortunately, that log-concavity fails for general
$\t>0$, even for the limit functions \tht{1.5}. Indeed, the coefficients
in the expansion
$$
\frac1{(1-z)^\t}=1+\t\,z+\frac{\t(\t+1)}2 z^2 + \dots 
$$
do not form a log-concave sequence if $\t<1$. In this paper 
(in the proof of the ``only if'' part of Theorem 1.1)
we shall use a different approach, where the main role is played by
the simple yet effective Lemma 5.2. 

\subhead
1.4.5~$S(\infty)$ and Kerov's conjecture (\cite{VK1,K,KOO})
\endsubhead

The work of A.~Vershik and S.~Kerov on the characters of $U(\infty)$
had very much in common with their earlier work \cite{VK1} on the
characters of $S(\infty)$. The problems considered in this
paper do have their ``symmetric group counterparts'', too.  
Those are very similar, but less technically involved. 
For example, the analog of the problem, addressed in Theorem 1.3, is
is the following problem: describe all homomorphisms
$$
\sigma: \L \to \C
$$
such that
$$
\forall\l\quad \sigma(P_\l(x;\t))\ge 0\,. \tag 1.23
$$
Here $\L$ denotes the $\C$-algebra of symmetric functions in
infinitely many variables. In \cite{K}, S.~Kerov conjectured an
explicit description of all homomorphisms $\L \to \C$ that are
positive on the Macdonald polynomials. In particular, it implies 
that all 
homomorphism $\sigma$ satisfying \tht{1.23} are given
 by extended symmetric functions
(see Section 2.7)
$$
\L\owns f @>\quad\sigma\quad>>  f(\a;\b;\g;\t) \in \C\,, 
$$
where $\a=(\a_1\ge\a_2\ge\dots\ge 0)$, $\b=(\b_1\ge\b_2\ge\dots\ge 0)$
and $\g\ge 0$ are arbitrary. This particular case of
Kerov's conjecture was proved in \cite{KOO} using
the approximation techniques
proposed in \cite{V,VK1} and the $P^*$-functions
machinery. That paper contains also a discussion
of other equivalent forms of this problem and a
symmetric group analog of Theorem 1.1. 
In the present paper we shall use the 
estimate established in Theorem 7.1 of \cite{KOO}.

\subhead
1.4.6~Other series of root systems
\endsubhead

In our next paper we shall prove analogs of the results obtained
here for other classical series of root systems. The necessary
$P^*$-functions techniques were developed in \cite{Ok10}, see also
our paper \cite{OO3}. Note also that the analogs of the Vershik-Kerov
theorem for the groups $O(\infty)$ and $Sp(\infty)$ were
obtained by R.~Boyer in \cite{B2}.

\subhead
1.5 Acknowledgments
\endsubhead

During the last years we had  fruitful discussions with
I.~Cherednik, R.~Howe, A.~A.~Kirillov, F.~Knop, A.~Lascoux, I.~Macdonald, 
A.~Molev, M.~Nazarov, S.~Sahi, and our 
colleagues from the Institute for Problems of Information
Transmission. We are very
grateful to all of them.
We are especially indebted to A.~Vershik and
S.~Kerov who very much influenced and inspired our work. 
S.~Kerov also read the preliminary version of this 
text and made a large number of valuable critical remarks.

The authors were supported by the Russian Basic Research Foundation
(grant 95-01-00814).
The first author enjoyed the hospitality of the Institute for Advanced
Study in Princeton and the Mathematical Sciences Research Institute
in Berkeley (which was made possible thanks to NSF grants 
DMS--9304580 and DMS--9022140 respectively).
 

\head
2.~Jack polynomials and shifted Jack polynomials
\endhead

In this section we fix notation and gather some 
fundamental results. We shall mostly follow the standard
Macdonald's notation, with a few exceptions. 
Most importantly, as in \cite{OO2}, our fixed 
positive parameter $\t>0$ is inverse to the Macdonald's 
parameter $\a=1/\t$.

\subhead 2.1~Orthogonality (\cite{M1,St})
\endsubhead 
We denote by $\L$ the $\C$-algebra of symmetric functions
(in infinitely many variables). The Jack polynomials $\Pl(x;\t)$
form an orthogonal basis of $\L$ with respect to the following
inner product:
$$
(p_\l,p_\mu):=\d_{\l,\mu} z_\l \,\t^{-\ell(\l)}\,, \tag 2.1 
$$
where $\l$ and $\mu$ are two partitions, $p_\l$ is the following
symmetric function
$$
p_\l:=\prod_{k\le\ell(\l)}
\left(\sum_i x_i^{\l_k}\right) \,,
$$
$\ell(\l)$ stands for
the number of parts of $\l$, and the factor $z_\l$ is defined in
\cite{M1}, Section I.2.
Let 
$$
Q_\l(x;\t):=\frac1{(\Pl,\Pl)} \Pl(x;\t) 
$$
be the dual (with respect to the inner product \tht{2.1}) basis of $\L$.
The orthogonality relation $(\Pl,\Qm)=\d_{\l,\mu}$ is equivalent to the following
Cauchy-type identity:
$$
\align
\prod_{i,j} (1- x_i y_j)^{-\t}&=\sum_\l \Pl(x) \, Q_\l(y) \,, 
\tag 2.2
\\
&=\sum_\l Q_\l(x) \, \Pl(y) \,, 
\endalign
$$

\subhead 2.2~Interpolation (\cite{S1,OO,Ok6,KS,Kn,S2,Ok8})
\endsubhead 

We denote by $\Lt(n)\subset\C[x_1,\dots,x_n]$ the subalgebra of
polynomials, symmetric in variables $x_i-\t i$, and denote by
$$
\Lt:=\varprojlim \Lt(n),
$$
the inverse limit of $\Lt(n)$ as filtered (by degree of polynomials)
algebras with respect to homomorphisms:
$$
\align
\Lt(n+1)\quad  &\to  \quad \Lt(n) \,, \\
f(x_1,\dots,x_{n+1})&\mapsto f(x_1,\dots,x_n,0)\,.
\endalign
$$
Note that, in particular, the value $f(\l)$ is well defined for
any partition $\l$ and any $f\in\Lt$. The {\it shifted} (or 
{\it interpolation}) Jack polynomials $\Pm(x;\t)$ is the 
unique polynomial satisfying the following Newton interpolation conditions:
\roster
\item $\Pm(x;\t)\in\Lt$\,,
\item $\deg\Pm(x;\t)=|\mu|$\,,
\item $\Pm(\l;\t)=0$ unless $\mu\subset\l$\,,
\item $\Pm(\mu;\t)=H(\mu)$\,,
\endroster
where $H(\mu)$ is a normalization constant defined as follows. Recall that
for a square  $s=(i,j)\in\mu$ in the diagram of a partition $\mu$ the numbers
$$
\alignat2
&a(s)=\mu_i-j,&\qquad &a'(s)=j-1,\\
&l(s)=\mu'_j-i,&\qquad &l'(s)=i-1,
\endalignat
$$
are called arm-length, arm-colength, leg-length, and
leg-colength, respectively. We set
$$
H(\mu):=\prod_{s\in\mu}(a(s)+\t\, l(s)+1) \,.
$$
Set also
$$
H'(\mu):=\prod_{s\in\mu}(a(s)+\t\, l(s)+\t) \,.
$$
Then we have 
$$
(\Pl,\Pl)=H(\l)/H'(\l)\,.  
$$

\subhead 2.3 Branching rules
\endsubhead
One has the following branching rule for the Jack polynomials 
(see \cite{St} and also \cite{M1}) 
$$
\Pl(x_1,\dots,x_n;\t)=\sum_{\mu\prec\l} \plm \, x_1^{|\l/\mu|} \, 
P_\mu(x_2,\dots,x_n;\t)\,, \tag 2.3
$$
where $\mu\prec\l$ stands for the inequalities  of interlacing 
$$
\l_1 \ge \mu_1 \ge \l_2 \ge \mu_2 \ge \dots
\ge \mu_{n-1} \ge \l_n\,,
$$
and (we use the standard notation $(t)_m=t(t+1)\cdots(t+m-1)$)
$$
\plm:=\prod_{1\le i \le j <n} 
\frac{(\mu_i-\mu_j+\t(j-i)+\t)_{\mu_j-\l_{j+1}}}
{(\mu_i-\mu_j+\t(j-i)+1)_{\mu_j-\l_{j+1}}}
\frac{(\l_i-\mu_j+\t(j-i)+1)_{\mu_j-\l_{j+1}}}
{(\l_i-\mu_j+\t(j-i)+\t)_{\mu_j-\l_{j+1}}}
$$
is a positive coefficient 
$$
\plm > 0\,, \quad \mu\prec\l\,,
$$ 
which is rational in $\t$ and equals 1 for
$\t=1$ (the Schur functions case).

The analog of \tht{2.3} for the interpolation Jack polynomials is
the following formula (see \cite{Ok8}):
$$
P^*_\l(x_1,x_2,\dots;\t)=\sum_{\mu\prec\l} \plm 
\left(\prod_{s\in\l/\mu}\left(x_1-a'(s)+\t\,l'(s)\right) \right)
\Pm(x_2,\dots;\t)\,. \tag 2.4 
$$
Here $\l$ is a partition whereas in \tht{2.3} it can be a
signature (see section 2.6). 
In particular, one has (see \cite{KS} and also \cite{Ok8}):
$$
P^*_\l(x;\t)=\Pl(x;\t)+\dots\,, \tag 2.5  
$$
where dots stand for the lower degree terms. 

Such a close relationship between the orthogonal polynomials $P_\mu$
and the interpolation polynomials $\Pm$ seems to be a rather 
remarkable phenomenon (notice that their definitions have nothing
in common). The interplay between $P_\mu$ and  $\Pm$ will play the central
role in this paper. 

\subhead 2.4 Binomial formula (\cite{OO2})
\endsubhead 
Given a partition $\mu$ and a number $t$ set
$$
(t)_\mu = \prod_{s\in\mu}(t+a'(s)-\t\, l'(s))\,.
$$
If $\mu=(m)$ then $(t)_\mu=(t)_m$ is the standard shifted factorial.
We have the following equivalent formulas 
$$
\align
\frac{P_\l(1+x_1,\dots,1+x_n;\t)}
{P_\l(1,\dots,1;\t)} &=
\sum_\mu 
\frac{ \Pm(\l;\t)\, Q_\mu(x_1,\dots,x_n;\t)}
{(n\t)_\mu} \tag 2.6 \\
&=
\sum_\mu 
\frac{ Q^*_\mu(\l;\t)\, P_\mu(x_1,\dots,x_n;\t)}
{(n\t)_\mu}\,,
\endalign 
$$
where $Q^*_\mu(x):=\Pm(x;\t)/(P_\mu,P_\mu)$. Here 
$\l$ can be a signature but $\mu$ is a partition. 
Recall that by $\P_\l$ we denote the normalized Jack 
rational function
$$
\P_{\l}(z_1,\dots,z_n;\t):=
\frac{P_{\l}(z_1,\dots,z_n;\t)}
{P_{\l}(1,\dots,1;\t)} \,.
$$
It follows that
$$ 
\P_\l(z_1,\dots,z_k,\overbrace{1,\dots,1}^{
\text{$n-k$ times}};\t)
=
\sum_{\ell(\mu)\le k} 
\frac{ Q^*_\mu(\l;\t)\, P_\mu(z_1-1,\dots,z_k-1;\t)}
{(n\t)_\mu}\,. \tag 2.7
$$
This expansion \tht{2.7} will play a fundamental role in this paper. 

\subhead 2.5 Generating functions 
\endsubhead 
By definition, set
$$
g_k(x;\t):=Q_{(k)}(x;\t)\,, \quad 
g^*_k(x;\t):=Q^*_{(k)}(x;\t)\,, \quad k=0,1,2,\dots\,.
$$
{F}rom the general formulas \tht{2.3-4} we have 
$$
\align
g_k(x;\t)&=
\sum_{1\le i_1 \le \dots \le i_k }
\frac{(\t)_{m_1}(\t)_{m_2}\cdots }{m_1!\,m_2!\, \cdots}
\,x_{i_1} x_{i_2} \cdots x_{i_k}\,, \\
g^*_k(x;\t)&=
\sum_{1\le i_1 \le \dots \le i_k }
\frac{(\t)_{m_1}(\t)_{m_2}\cdots }{m_1!\,m_2!\, \cdots}
\,(x_{i_1}-k+1) \cdots (x_{i_{k-1}}-1) x_{i_k}\,, \tag 2.8
\endalign
$$
where $m_l:=\#\{r\,|\, i_r=l\}$ stand for the multiplicities with
which the numbers  $l=1,2\dots$ occur in $i_1,\dots,i_k$. 

Consider the following generating functions:
$$
G(x;t):=\sum_{k\ge 0} g_k(x)\, t^k\,, \quad
G^*(x;u):=\sum_{k\ge 0} \frac{g^*_k(x)}{u(u-1)\cdots(u-k+1)}\,.
$$
Observe that the last series is a well defined element of 
$\Lt[[u^{-1}]]$ and
$$
G(x;t)=\lim_{a\to\infty} G^*(ax;a/t) \,.
$$
It follows from the Cauchy identity \tht{2.2} that 
$$
G(x;t)=\prod_i (1-t x_i)^{-\t}\,. \tag 2.9 
$$
The interpolation analog of \tht{2.9} is the following evaluation
$$
G^*(x;u)=\prod_i 
\frac{\G(x_i-u-\t\,i)}{\G(x_i-u-\t\,i+\t)}
\frac{\G(-u-\t\,i+\t)}{\G(-u-\t\,i)} \,. \tag 2.10
$$
Originally, we obtained the formula \tht{2.10} as an
auxiliary statement in the proof of the formula \tht{2.8}
for the polynomials $g^*_k(x;\t)$ (unpublished). That proof
of \tht{2.10} was based on iterated use of the Gauss summation formula;
it is reproduced in the Appendix to this paper (see Section 7).
Alternatively, the formula \tht{2.10} can be obtained by taking the
limit $q\to 1$ in its $q$-analog proved in \cite{Ok9} or by
using the very same argument as used in the proof of the
generating functions \tht{12.3} in \cite{OO} and \tht{2.9} in \cite{Ok9}. Also,
it can be regarded as a particular case of the binomial formula
for $P^*$-functions, see \cite{Ok9}. 

\subhead 2.6 Partitions and signatures  
\endsubhead 
We call any non-increasing sequence of integers 
$$
\l=(\l_1\ge \dots \ge \l_n)\,, \quad n=1,2,\dots 
$$
a {\it signature}. The number $n$ is called the {\it
length} of the signature $\l$. Using the relation
$$
P_{(\l_1+1,\dots,\l_n+1)} (x_1,\dots,x_n;\t)= 
\left(
\prod_{i=1}^n x_i 
\right)\,P_\l(x;\t)
$$
one can define rational Jack functions $\Pl(x;\t)$ for any
signature $\l$. These functions enjoy the same branching rules
and binomial expansions as the Jack polynomials. It is important to
notice that the definition of the rational Jack polynomials 
makes sense only for finitely many variables (namely, for $n$ 
variables where $n$ is the length of the signature $\l$).

To each signature $\l$ one can associate two partitions $\lp$
and $\lm$ as follows. Suppose $\l_p\ge 0\ge \l_{p+1}$ for some
$p=1,\dots,n$. Then, by definition,
$$
\align
\lp&=(\l_1,\dots,\l_p)\,,\\
\lm&=(-\l_n,\dots,-\l_{p+1})\,.
\endalign
$$
This decomposition has the following property
$$
G^*(\l;u)=G^*(\lp;u) \, G^*(\lm;-u-\t\,n-1) \,, \tag 2.11
$$
which follows from \tht{2.10} and the following elementary identity
$$
\multline
\frac{\G(x-u-\t)}{\G(x-u)}
\frac{\G(-u)}{\G(-u-\t)}=
\frac{(x-u)\cdots(-1-u)}
{(x-u-\t)\cdots(-1-u-\t)}=\\
=\frac{(u+1)\cdots(u-x)}
{(u+\t+1)\cdots(u+\t-x)}=
\frac{\G(-x+u+1)}{\G(-x+u+\t+1)}
\frac{\G(u+\t+1)}{\G(u+1)}\,, \quad x\in{\Bbb Z_{\le 0}} \,. 
\endmultline
$$
A similar argument works for $x\in{\Bbb Z_{\ge 0}}$. 
Observe, however, that the above identity is false if $x\notin{\Bbb Z}$.

\subhead 2.7 Extended symmetric functions  
\endsubhead 
Let $f$ be an element of $\L$. Recall that we have fixed
a positive real parameter $\t>0$.  

Given two sequence of variables  $\a=(\a_1,\a_2,\dots)$,
$\b=(\b_1,\b_2,\dots)$ and one extra variable $\g$ 
we define (cf.\ \cite{K}) the {\it extended} symmetric function 
$f(\a;\b;\g;\t)$ as the result of specialization 
$$
\L\to\C[\a,\b,\g]
$$
given by
$$
\align
p_1 &\mapsto \sum \a_i + \sum \b_i + \g\,, \\
p_k &\mapsto \sum \a_i^k + (-\t)^{k-1} \sum \b_i^k\,, \quad k\ge 2 \,.
\endalign
$$
Equivalently, this specialization can be described by its action
on the generating series $G(t)$:
$$
G(t)\mapsto e^{\g t \t} \prod_i \frac{(1+t \t \, \b_i)}
{(1-t \a_i)^{\t}} \,.
$$
We shall use this definition in the situation when 
$\a_i,\b_i,\g\in\R$ and the series 
$$
\sum \a_i, \,\, \sum \b_i < \infty
$$
converge absolutely. 
In the sequel we write simply $f(\a;\b;\g)$ instead of 
$f(\a;\b;\g;\t)$

More generally, given 4 sequences $\ap,\am,\bp,\bm$ of variables 
and two extra variables  $\gp,\gm$ we define the 
{\it doubly extended} symmetric function 
$$
f\ABG
$$
as the result of specialization defined by
$$
G(t)\mapsto e^{\gp \t t + \gm \t t'} \prod_i 
\frac{(1+t\t \bpi)}
{(1-t \api)^{\t}} 
\frac{(1+t'\t \bmi)}
{(1-t' \ami)^{\t}}\,, \tag 2.12 
$$
where $t':= -t/(1+\t t)$, or, equivalently 
$$
(1+\t t)(1+\t t')=1\,.
$$

\head
3.~Asymptotic properties of Vershik-Kerov sequences of 
signatures
\endhead

In this section we prove the following 
\proclaim{Theorem 3.1} Suppose $\l(n)$ is a VK sequence of
signatures with parameters $\a^\pm_i$, $\b^\pm_i$, and $\g^\pm$.
Suppose $f^*\in\Lt$ and let $f\in\L$ be the highest degree 
term of $f^*$. Then
$$
\lim_{n\to\infty} \frac{\tsize f^*(\l(n))}{\tsize n^{\deg f^*}} = 
f\ABG\,. \tag 3.1
$$
In particular, the limit  depends on the highest degree term $f$ of $f^*$
only.
\endproclaim 

\demo{Proof} In the particular case when $\l(n)$ is, in fact,
a sequence of partitions this theorem was established in \cite{KOO}.
Namely, it follows from the estimate stated in Theorem 7.1 of 
\cite{KOO} and Lemma 5.2 of \cite{KOO}. Notice that in \cite{KOO}
a different convention about the letter $n$ was used: $n$ stood
there for the number of squares $|\l(n)|$ in the partition $\l(n)$. Here
we have weaker assumptions, namely
$$
\ell(\l(n))\le n\,,\quad\text{and}\quad\exists\lim_{n\to\infty}
\frac{|\l(n)|}n\,.
$$
It is easy to see that this difference in assumptions is unessential. 

We shall deduce
the general statement \tht{3.1} from that particular case. 
Observe that it suffices to prove \tht{3.1} for
$$
f^*=g^*_k\,,\quad  \quad f=g_k\,, \quad k=1,2,\dots\,,
$$
because the functions $\{g_k\}$ are homogeneous 
generators  of the algebra $\L$.
To simplify notation, set
$$
g(k,n)=g^*_k(\l(n))\,, \quad g^\pm(k,n)=g^*_k(\l^\pm(n)) \,.
$$
Since $\l^\pm(n)$ are two VK sequences of partitions, we conclude that
the following limits exist:
$$
\align
g^+(k)&:=\lim_{n\to\infty} \frac{g^+(k,n)}{n^k} =
g_k(\ap;\bp;\gp)\,, \\
g^-(k)&:=\lim_{n\to\infty} \frac{g^-(k,n)}{n^k} =
g_k(\am;\bm;\gm)
\endalign
$$
and satisfy the equality
$$
\sum_{k} g^\pm(k) \,t^k = 
e^{\g^\pm \t t} \prod_i 
\frac{(1+t\t \b^\pm_i)}
{(1-t \a^\pm_i )^{\t}} \,. \tag 3.2 
$$

The identity \tht{2.11} can be rewritten as
$$
\multline
\sum_{k\ge 0} \frac{g(k,n)}{u(u-1)\cdots(u-k+1)} = \\
\left( \sum_{p\ge 0} \frac{g^+(p,n)}{u(u-1)\cdots(u-p+1)} \right)
\left( \sum_{q\ge 0} \frac{g^-(q,n)}{u'(u'-1)\cdots(u'-q+1)}
\right)\,,
\endmultline \tag 3.3  
$$
where
$$
u':=-u-\t\,n-1\,.
$$
We introduce a new variable $v$
$$
u=n v \,.
$$
Then \tht{3.3} becomes
$$
\multline
\sum_{k\ge 0} \frac{g(k,n)}{n^k}
\frac1{v(v-\frac1{n})\cdots(v-\frac{k-1}{n})} = \\
\left( \sum_{p\ge 0} \frac{g^+(p,n)}{n^p}
\frac1{v(v-\frac1{n})\cdots(v-\frac{p-1}{n})} \right)
\left( \sum_{q\ge 0} \frac{g^-(q,n)}{n^q} 
\frac1{(v'-\frac1{n})\cdots(v'-\frac{q}{n})}
\right)\,,
\endmultline \tag 3.4 
$$
where
$$
v':=-v-\t\,.
$$
Recall that we consider \tht{3.4} as an identity in  the formal
power series algebra $\C[[1/v]]$. We equip $\C[[1/v]]$ with
the topology of the coefficient-wise convergence. Clearly,
in  this topology
$$
\frac1{v-a_n} \to \frac1{v-a_\infty}\,, \quad n\to\infty\,,
$$
provided $a_n\to\a_\infty$. Therefore the
RHS of \tht{3.4} converges  to the following element of $\C[[1/v]]$
$$
\left( \sum_{p\ge 0} \frac{g^+(p)}{v^p} \right)
\left( \sum_{q\ge 0} \frac{g^-(q)}{(-v-\t)^q}
\right)\,.
$$
Comparing the like powers of $1/v$ in \tht{3.4} we conclude 
by induction on $k$ that the following limits exist
$$
g(k):=\lim_{n\to\infty} \frac{g(k,n)}{n^k}
$$
and satisfy the equality
$$
\sum_k \frac{g(k)}{v^k} = \left( \sum_{p\ge 0} \frac{g^+(p)}{v^p} \right)
\left( \sum_{q\ge 0} \frac{g^-(q)}{(-v-\t)^q}
\right)\,.
$$
Set $t=1/v$. Then we obtain from \tht{3.2} 
$$
\sum_k g(k)\, t^k 
=e^{\gp \t t + \gm \t t'} \prod_i 
\frac{(1+t\t \bpi)}
{(1-t \api)^{\t}} 
\frac{(1+t'\t \bmi)}
{(1-t' \ami)^{\t}}\,,
$$
where 
$$
t'=-\frac1{v+\t}= -\frac{t}{1+\t t}\,.
$$ 
This concludes the proof of the theorem. \qed 
\enddemo

\head
4.~Sufficient conditions of regularity
\endhead

In this section we prove the following
\proclaim{Theorem 4.1}
Suppose $\l(n)$ is a VK sequence of signatures with parameters
$\a^\pm_i$, $\b^\pm_i$, and $\g^\pm$. Then $\l(n)$ is
regular and infinitesimally regular. Set
$$
\P_{\a,\b,\g}(z_1,z_2,\dots) := \phi_{\a,\b,\g}(z_1) 
\phi_{\a,\b,\g}(z_2) \dots\,,
$$
where the function $\phi_{\a,\b,\g}(z)$ is the following product 
$$
\phi_{\a,\b,\g}(z):=e^{\gp (z-1) + \gm (z^{-1}-1)} \prod_i 
\frac{(1+ \bpi (z-1))}
{(1- \api (z-1)/ \t)^{\t}} 
\frac{(1+\bmi (z^{-1}-1))}
{(1-\ami  (z^{-1}-1)/ \t)^{\t}}\,. 
$$
Then
$$
\P_{\l(n)}(z_1,\dots,z_k,1,\dots ;\t)
\to \P_{\a,\b,\g}(z_1,\dots,z_k,1,\dots)
$$
uniformly on each torus $\T^k$, $k=1,2,\dots$.
\endproclaim

We shall deduce the above theorem from Theorem 3.1 of the previous
section and the following
\proclaim{Lemma 4.2}
Fix some $k=1,2,\dots$ and suppose that $\psi_1,\psi_2,\dots,\psin$
are smooth functions on the torus $\T^k$ such that
\roster
\ite a $\psi_n(1,\dots,1)=1\,, \quad n=1,2,\dots,\infty$\,;
\ite b the functions $\psi_n$ are positive definite; 
\ite c $\psin$ is a real-analytic function on $\T^k$;
\ite d the Taylor expansions of $\psi_n$ about the
point $(1,\dots,1)\in\T^k$ converge coeffi\-ci\-ent-wise to 
the Taylor expansion of $\psin$;
\endroster
then $\psi_n\to\psin$ uniformly on $\T^k$. 
\endproclaim 

\demo{Proof of Theorem 4.1} We have to verify that the functions
$$
\align
\psi_n(z_1,\dots,z_k):&=\P_{\l(n)}(z_1,\dots,z_k,
\underbrace{1,\dots,1}_{\text{$n-k$ times}};\t)\,, \\
\psin(z_1,\dots,z_k):&=
\P_{\a,\b,\g}(z_1,\dots,z_k,1,\dots)
\endalign 
$$
satisfy all conditions of Lemma 4.2. 
Obviously, $\psi_n$ are smooth and 
$\psi_n(1,\dots,1)=1$ for all $n$. 

Since the coefficients $\plm$ in the branching rules \tht{2.3}
are nonnegative, it follows that $\psi_n$ are polynomials
in $z^{\pm}_i$ with nonnegative coefficients and hence
positive definite functions. 

Since the functions $\psi_n$ are symmetric and the Jack 
polynomials form a linear basis of $\L$, we can rewrite
the Taylor expansion of $\psi_n$ as an expansion in
polynomials 
$$
P_\mu(t_1,\dots,t_k;\t)\,, \quad t_i=z_i-1\,,\quad \ell(\mu)\le k \,.
$$
This expansion is given by the formula \tht{2.7}, which in our
current notation reads
$$
\psi_n(z_1,\dots,z_k)=\sum_{\ell(\mu)\le k} 
\frac{Q^*_\mu(\l(n))}{(\t\,n)_\mu}  P_\mu(t_1,\dots,t_k;\t)\,. 
\tag 4.1  
$$
It is clear that
$$
\frac{Q^*_\mu(\l(n))}{(\t\,n)_\mu} \sim
{\t^{-|\mu|}}\,\frac{Q^*_\mu(\l(n))}{n^{|\mu|}}\,,
\quad n\to\infty\,.
$$
By assumption the sequence $\l(n)$ is VK; therefore, by \tht{2.5} and 
Theorem 3.1  we have
$$
\frac{Q^*_\mu(\l(n))}{(\t\,n)_\mu} \to {\t^{-|\mu|}}\, Q_\mu\ABG\,, \quad 
n\to\infty \,.
$$

Now it suffices to verify that
$$
\psin(z_1,\dots,z_k)
=\sum_{\ell(\mu)\le k} 
{\t^{-|\mu|}}\, Q_\mu\ABG  P_\mu(t_1,\dots,t_k;\t)\,.
$$
But this follows from the Cauchy identity \tht{2.2} which 
implies
$$
\multline
\sum_{\ell(\mu)\le k} 
{\t^{-|\mu|}}\, Q_\mu(x_1,x_2,\dots)  P_\mu(t_1,\dots,t_k;\t) 
=
\\
=\prod_{i=1}^\infty \prod_{j=1}^k \frac1{(1-x_i t_j/\t)^\t}
=\prod_{j=1}^k G(x_1,x_2,\dots;t_j/\t) 
\endmultline 
$$
and from the definition \tht{2.12} of the specialization 
$$
Q_\mu\ABG \,.
$$
This concludes the proof of the theorem. \qed
\enddemo 
  
We have the following corollary of the above proof:
\proclaim{Corollary 4.3} We have the following equivalence 
between two conditions listed in Theorem 1.2 
$$
\text{\tri{ii}}\Longleftrightarrow
\text{\tri{iv}}\,.
$$
\endproclaim 
\demo{Proof} Follows from \tht{4.1} and the fact 
that the polynomials $Q_\mu$ form a homogeneous linear
basis of the algebra $\L$.\qed
\enddemo

Now it remains to establish Lemma 4.2. 
The proof of the lemma is quite standard. To avoid multiindices we assume
$k=1$; the argument in the general case is the same.

\demo{Proof of Lemma 4.2}
By \tri{a} and \tri{b} we have
\footnote
{
Here and below it will be convenient for us to write sums as
integrals with respect to a discrete  measure. Most of our
computations (such as Lemma 5.2  below) apply to arbitrary measures
on $\R$.
}   
$$
\psi_n(\e s)=\int \e{s\xi}\, M_n(d\xi)\,, \quad s\in\R/\Z\,,
\quad \xi\in\Z\,,
$$
where $M_n$ is a certain probability measure supported on 
$\Z=\T^\wedge$. By assumption, all functions $\psi_n$ are
smooth and hence all moments of all
measures $M_n$ are finite.
By \tri{d} the limits
$$
m_l:=\lim \int \xi^l \, M_n(d\xi)\,, \quad n\to\infty\,,
$$
exist. In particular, we have certain upper bounds
$$
\int \xi^l \, M_n(d\xi) \le c(l)
$$
uniformly in $n$. The following Chebyshev-type uniform tail estimates
$$
\int_{|\xi|>a} \xi^{2l}\, M_n(d\xi) \le
\int_{|\xi|>a} \frac{\xi^{2l+2}}{a^2}\, M_n(d\xi) \le
\frac{c(2l+2)}{a^2}\,,
$$
where $a$ is arbitrary positive, 
imply that the family $\{M_n\}$ is tight (and hence
relatively compact in the weak convergence topology)
and, moreover,
$$
\int \xi^l \, M_\infty(d\xi) = m_l 
$$
for any limit point $M_\infty$ of the set $\{M_n\}$. 
Note that $M_\infty$ is a probability measure.

By \tri{c} the series 
$$
\psi_\infty(\e s)=\sum_{l\ge 0} m_l \frac{(2\pi i s)^l}{l!}
$$
converges in some neighborhood of $s=0$. Hence (see e.g.\ 
Theorem II.12.7 in \cite{Sh})  the moments $\{m_l\}$ determine
the measure $M_\infty$ uniquely and (see the proof of the 
aforementioned theorem) the Fourier transform of $M_\infty$
is real-analytic and hence equals $\psi_\infty$ everywhere. 

It follows that the set $\{M_n\}$ has a unique limit point
$M_\infty$, whence
$$
M_n @>\text{\,weak\,}>> M_\infty\,, \quad n\to\infty
$$
and hence
$$
\psi_n \to \int \e{s\xi}\, M_\infty(d\xi) = \psi_\infty
$$
uniformly on compact sets. Since $\T$ is compact, the lemma
follows. \qed
\enddemo

\head
5.~Necessary conditions of regularity
\endhead

\subhead 
5.1 The ``only if'' part of Theorem 1.1
\endsubhead

Given a vector $\l=(\l_1\ge\dots\ge\l_n)\in\R^n$ set,
by definition, 
$$
\N{\l}^2:= \sum \l_i^2 + \t\left(\sum \l_i\right)^2 + \t(\l,2\r)\,,
$$
where $2\r:=(n-1,n-3,\dots,3-n,1-n)$ and hence
$$
(\l,2\r)=\sum_{i<j} (\l_i-\l_j) \ge 0 \,.
$$
The function $\N{\l}$ will play the role of a ``norm'' of $\l$.

The main theorem which we shall prove in this section is the 
following

\proclaim{Theorem 5.1} If $\l(n)$ is a sequence of signatures and 
the following limit  
$$
\lim_{n\to\infty} \P_{\l(n)} (z,\underbrace{1,\dots,1}_{\text{$n-1$ times}})\,, 
\quad z\in\T\,, \tag 5.1
$$
exists point-wise and is continuous at $z=1$ then
the sequence $\l(n)$ is a Vershik-Kerov sequence. 
\endproclaim 

In particular, the hypothesis of Theorem 5.1 is clearly satisfied
if the sequence $\l(n)$ is regular. Therefore, the combination of
this theorem with Theorem 4.1 implies Theorem 1.1.

We shall prove Theorem 5.1 in three
steps: we shall prove that its hypothesis implies that
\roster
\item $\N{\l(n)}=O(n)\,,\quad  n\to\infty$,
\item $|\l^\pm(n)|=O(n)\,,\quad  n\to\infty$, \quad  and, finally, 
\item $\l(n)$ is a Vershik-Kerov sequence.
\endroster
To justify the first step  we shall need the following 

\proclaim{Lemma 5.2} Let $\{M_n\}$ be a tight family of probability
measures on $\R$ with finite $4$-th moments
$$
\int x^4 \,M_n(dx) < \infty \,, 
$$
and suppose that
$$
\int x^2 \,M_n(dx) \to \infty\,.
$$
Then the following ratio grows to infinity 
$$
{\dsize\int x^4 \,M_n(dx)}\left/{\dsize\left(\int x^2 \,M_n(dx)
\right)^2}\right.\to\infty\,.
$$
\endproclaim
\demo{Proof} Fix some $\ep>0$ and show that
$$
4\ep \int x^4 \,M_n(dx) \ge \left(\int x^2 \,M_n(dx)
\right)^2\,, \tag 5.2 
$$
provided $n$ is sufficiently large. Since $\{M_n\}$ is tight we can
find (and fix) some $a>0$ such that
$$
\int_{|x|\ge a} M_n(dx) \le \ep
$$
for all $n$. By Cauchy inequality we have
$$
\left(\int_{|x|\ge a} x^2 \,M_n(dx)
\right)^2 \le \ep \int_{|x|\ge a}x^4\, M_n(dx)
\le \ep \int x^4 \,M_n(dx)\,. \tag 5.3 
$$
On the other hand 
$$
\int x^2 \,M_n(dx) \le a^2 + 
\int_{|x|\ge a} x^2 \,M_n(dx)
$$
and since $\int x^2 \,M_n(dx)\to\infty$ we have
$$
\int_{|x|\ge a} x^2 \,M_n(dx) \ge \frac12 \int x^2 \,M_n(dx) \tag 5.4 
$$
for all sufficiently large $n$. Substituting \tht{5.4} into \tht{5.3}
we obtain \tht{5.2}. \qed
\enddemo

Now we begin with the proof of the Theorem 5.1.

\demo{Proof of Theorem 5.1}

\noindent{\it Step 1}. Our goal now is to prove that
$$
\N{\l(n)}=O(n)\,,\quad  n\to\infty\,. \tag 5.5 
$$
We set
$$
\phi_n(z):=\P_{\l(n)} (z,\underbrace{1,\dots,1}_{\text{$n-1$ times}})
$$
and denote by $\pin$ the limit function \tht{5.1}. 

As in the previous section
we shall represent the function $\phi_n(z)$ as the Fourier
transform of some probability measure $M_n$ on $\Z$
$$
\phi_n(z)= 
\int z^\xi M_n(d\xi)\,, \quad z\in\T^1\,,
\quad \xi\in\Z\,.
$$
By the continuity property of characteristic functions 
(see, for example, \cite{Sh}, Theorem III.3.1) we have
$$
M_n @>\text{\quad weak\quad}>> M_\infty\,,
$$
where 
$$
\pin(z)=\int z^\xi M_\infty(d\xi)\,.
$$

To simplify notation we shall write $\l$ in place of $\l(n)$. 
The binomial expansion
$$
\phi_n(z)=\sum_{k\ge 0} \frac{g^*_k(\l)}{(n\t)_k} \, (z-1)^k
$$
implies that 
$$
\int \xi(\xi-1)\dots(\xi-k+1)\, M_n(d\xi)= k!\,
\frac{g^*_k(\l)}{(n\t)_k} \,,
$$
whence the $k$-th moment of $M_n$ is a linear combination
of the numbers
$$
\frac{g^*_1(\l)}{(n\t)_1},\quad\dots\quad,
\frac{g^*_k(\l)}{(n\t)_k}\,.
$$
For example, we have
$$
\int \xi^2\, M_n(d\xi) = 
\frac{g^*_1(\l)}{n\t}+ 2\frac{g^*_2(\l)}{n\t(n\t+1)}\,.
$$
The formula \tht{2.8} specializes to 
$$
\align
g^*_1(\l)&=\l_1+\dots+\l_n\,, \\
g^*_2(\l)&=
\t^2\sum_{i<j}(\l_i-1)\l_j+\frac{\t(1+\t)}2\sum_i (\l_i-1)\l_i \,.
\endalign
$$
After some simple algebra, we obtain
$$
\int \xi^2\, M_n(d\xi) = \frac{\N{\l}^2}{n(n\t+1)}\,. 
$$ 
Therefore, to complete the first step of the proof we have to show
that the $2$-nd moments of $M_n$ remain bounded as $n\to\infty$ 
$$
\int \xi^2\, M_n(d\xi) = O(1)\,, \quad n\to\infty\,. \tag 5.6 
$$
Below in Theorem 5.4 we shall prove a general growth estimate on 
functions $|f(\l)|$, where $f\in\Lt$ is arbitrary.
It implies, in particular, that
$$
\frac{g^*_k(\l)}{n^k} = O\left(1+\frac{\N{\l}^k}{n^k}\right)\,. 
$$
This yields that 
$$
\int \xi^k \, M_n(d\xi) 
= O\left(1+\frac{\N{\l}^k}{n^k}\right)\,,
$$
which together with Lemma 5.2 implies \tht{5.6} and, hence, \tht{5.5}.

\medskip
\noindent{\it Step 2}. Our goal now is to prove that
$$
\big(\N{\l(n)}=O(n)\big)
\Longrightarrow
\big(|\l^\pm(n)|=O(n)\big)\,. 
$$
Again, to simplify notation we shall write simply $\l$ in place of $\l(n)$. 
By definition of $\N{\l}$ we have 
$$
\align
\sum \l_i &= O(n)\,, \tag 5.7\\
\sum_{i<j} (\l_i-\l_j) &= O(n^2) \,. \tag 5.8 
\endalign
$$
We can find two integers $p=p(n)$ and $q=q(n)$ such that
$p+q=n$, $\ell(\lp)\le p$, and $\ell(\lm)\le q$. 
Then we have
$$
\multline
\sum_{i<j} (\l_i-\l_j) \ge 
\sum_{i\le p<j} (\l_i-\l_j)=\\
= q|\lp|+p|\lm|= 
\frac n2(|\lp|+|\lm|) + \frac{q-p}2\sum \l_i\,,
\endmultline  \tag 5.9 
$$
where we used the equalities 
 $\sum\l_i=|\lp|-|\lm|$ and $p+q=n$\,. 
Since $q-p=O(n)$ we obtain from \tht{5.9} and \tht{5.7-8}
$$
n (|\lp|+|\lm|) = O(n^2) \,,
$$
which is equivalent to $|\l^\pm(n)|=O(n)$.

\medskip
\noindent{\it Step 3}. Now we finish the proof of the theorem. 
The standard compactness argument yields that any sequence $\{\l(n)\}$
of signatures such that $|\l^\pm(n)|=O(n)$ contains a VK subsequence.
Therefore, in order to show that our sequence $\{\l(n)\}$ is VK
it suffices to show that if 
$$
\{\dd{\l}(n)\},  \{\ddd{\l}(n)\} \subset \{\l(n)\} \tag 5.10 
$$
are two VK subsequences then the corresponding VK parameters are equal:
$$
(\dd{\a}^+,\dd{\a}^-,\dd{\b}^+,\dd{\b}^-,\dd{\g}^+,\dd{\g}^-)=
(\ddd{\a}^+,\ddd{\a}^-,\ddd{\b}^+,\ddd{\b}^-,\ddd{\g}^+,\ddd{\g}^-)
\,. \tag 5.11
$$
Take two arbitrary VK subsequences \tht{5.10}.
By Theorem 4.1 the corresponding subsequences $\{\dd{\phi}_n\}$ and 
$\{\ddd{\phi}_n\}$ converge uniformly and to the same limit:
$$
\multline 
e^{\dd{\g}^+ (z-1) + \dd{\g}^- (z^{-1}-1)} \prod_i 
\frac{(1+ \dd{\b}^+_i (z-1))}
{(1- \dd{\a}^+_i (z-1)/ \t)^{\t}} 
\frac{(1+\dd{\b}^-_i (z^{-1}-1))}
{(1-\dd{\a}^-_i  (z^{-1}-1)/ \t)^{\t}}=\\
=
e^{\ddd{\g}^+ (z-1) + \ddd{\g}^- (z^{-1}-1)} \prod_i 
\frac{(1+ \ddd{\b}^+_i (z-1))}
{(1- \ddd{\a}^+_i (z-1)/ \t)^{\t}} 
\frac{(1+\ddd{\b}^-_i (z^{-1}-1))}
{(1-\ddd{\a}^-_i  (z^{-1}-1)/ \t)^{\t}}\,. 
\endmultline \tag 5.12 
$$
Consider the singularities  of the two equal analytic functions in
\tht{5.12}. The LHS is regular in the annulus
$$
\frac{\dd{\a}^-_1}{\t+\dd{\a}^-_1}
< |z| < \frac{\t+\dd{\a}^+_1}{\dd{\a}^+_1}
$$
and becomes singular as $z\to \frac{\dd{\a}^-_1}{\t+\dd{\a}^-_1},
\frac{\t+\dd{\a}^+_1}{\dd{\a}^+_1}$. It follows that 
$$
\dd{\a}^+_1=\ddd{\a}^+_1\,, \quad 
\dd{\a}^-_1=\ddd{\a}^-_1\,.
$$
Thus, both sides of \tht{5.12} have common factors which can be
cancelled. Iterating the argument we obtain
$$
\dd{\a}^+=\ddd{\a}^+\,, \quad \dd{\a}^-=\ddd{\a}^-\,.
$$
Next, set 
$$
f(x):=1 -\frac{1}x\,.
$$ 
This function is increasing for $x>0$.
Consider the zeros of \tht{5.12}. This yields the 
equality of the  the following 
multisets (recall that 
$\bp_1+\bm_1\le 1$)
$$
\multline
\left\{
\dots\ge f(1-\dd{\b}^-_2) \ge f(1-\dd{\b}^-_1) \ge
f(\dd{\b}^+_1) \ge f(\dd{\b}^+_2) \ge \dots 
\right\}_{\dd{\b}^\pm_i\ne 0} = \\
\left\{
\dots\ge f(1-\ddd{\b}^-_2) \ge f(1-\ddd{\b}^-_1) \ge
f(\ddd{\b}^+_1) \ge f(\ddd{\b}^+_2) \ge \dots 
\right\}_{\ddd{\b}^\pm_i\ne 0} \,,
\endmultline
$$
whence
$$
\multline
\left\{
\dots\ge 1-\dd{\b}^-_2 \ge 1-\dd{\b}^-_1 \ge
\dd{\b}^+_1 \ge \dd{\b}^+_2 \ge \dots \right\} = \\
\left\{
\dots\ge 1-\ddd{\b}^-_2 \ge 1-\ddd{\b}^-_1 \ge
\ddd{\b}^+_1 \ge \ddd{\b}^+_2 \ge \dots \right\} \,,
\endmultline \tag 5.13  
$$
or, in words, the two sequences in \tht{5.13} coincide up to a possible
shift of indices. Let $r\in\Z$ be the shift of indices in \tht{5.13}.
Then
$$
e^{\dd{\g}^+ (z-1) + \dd{\g}^- (z^{-1}-1)} =
z^r e^{\ddd{\g}^+ (z-1) + \ddd{\g}^- (z^{-1}-1)}\,.
$$
It follows that $\dd{\g}^\pm=\ddd{\g}^\pm$ and $r=0$. Hence
$$
\dd{\b}^+=\ddd{\b}^+\,, \quad \dd{\b}^-=\ddd{\b}^-\,.
$$
This establishes \tht{5.11} and, thus, concludes the proof of the theorem
(assuming the truth of Theorem 5.4 to be established in the next
section). \qed
\enddemo
  
We have the following corollary of the above proof:
\proclaim{Corollary 5.3} We have the following implication 
between two conditions listed in Theorem 1.2 
$$
\text{\tri{ii}}\Longrightarrow
\text{\tri{iii}}\,.
$$
\endproclaim 
\demo{Proof} Argue as in the above theorem but without using
Lemma 5.2 and Theorem 5.4. \qed
\enddemo 

Now we have accumulated enough knowledge to prove Theorems 1.1 and 1.2
\demo{Proof of Theorem 1.1} Follows immediately from Theorems 4.1 and
5.1. \qed
\enddemo 
\demo{Proof of Theorem 1.2} We merely list all already established
implications
$$
\alignat2
\text{\tri{iii}}&\Longrightarrow
\text{\tri{i}, \tri{ii}}\,,& \qquad&\text{by Theorem 4.1}\,,\\
\text{\tri{ii}}&\Longleftrightarrow
\text{\tri{iv}}\,,& \qquad&\text{by Corollary 4.3}\,,\\
\text{\tri{i}}&\Longrightarrow
\text{\tri{iii}}\,,& \qquad&\text{by Theorem 5.1}\,,\\
\text{\tri{ii}}&\Longrightarrow
\text{\tri{iii}}\,,& \qquad&\text{by Corollary 5.3}\,,
\endalignat
$$
which prove that
$$
\text{\tri{i}}\Longleftrightarrow
\text{\tri{ii}}\Longleftrightarrow
\text{\tri{iii}}\Longleftrightarrow
\text{\tri{iv}}\,.
$$
Now the hypothesis of Theorem 5.1 is clearly weaker then \tri{i}. But
by Theorems 5.1 and 4.1 it also implies \tri{i}. This concludes the
proof. \qed
\enddemo

\subhead 
5.2 A growth estimate for $|f(\l)|$, $f\in\Lt$
\endsubhead

In this subsection we shall prove the following estimate
which we used in the proof of Theorem 5.1 (namely, in Step 2)
\proclaim{Theorem 5.4} For any $f\in\Lt$ there exists a constant $C_f$
such that
$$
\left|f(\l)\right|
\le C_f\,\max\left(\N{\l},n\right)^{\deg f} \tag 5.14 
$$
for all $n=1,2,\dots$ and all $\l=(\l_1\ge\dots\ge\l_n)\in\R^n$.
\endproclaim  

Although rough, this estimate is sufficient for our purposes.
First,  we establish the following
\proclaim{Lemma 5.5} Suppose $\l=(\l_1\ge\dots\ge\l_n)\in\R^n$ then
$$
\left|\sum_{i=1}^n \l_i (i-1)^m \right|
\le \left|\sum \l_i\right|\, n^m + 2 (\l,2\r) \, n^{m-1} \tag 5.15
$$
for all $m=0,1,2,\dots$. 
\endproclaim
\demo{Proof of Lemma 5.5}
Denote the linear function to be estimated by 
$$
l_m(\l):= \sum \l_i(i-1)^m\,.
$$
It is easy to see that using affine
transformations
$$
\l_i\mapsto a\l_i + b\,, \quad i=1,\dots,n\,,
$$
where $a,b\in\R$ and $a\ne 0$, one reduces the estimate \tht{5.15}
to the following estimate
$$
\left(
{\matrix
\l_1\ge\dots\ge\l_n\,,\\ \vspace{2 pt}
\sum \l_i = 0\,,\\ \vspace{2 pt}
(\l,2\r)=n
\endmatrix}
\right)
\Longrightarrow
|l_m(\l)|\le 2n^m
\,, \tag 5.16 
$$
which shall now be established. 
The vectors $\l$ satisfying the inequalities in the LHS of
\tht{5.16} form an $(n-2)$-dimensional simplex with extreme points
$$
\l^{(k)}=\Big(\underbrace{\frac 1k,\dots,\frac 1k
}_{\text{$k$ times}}, \underbrace{-\frac 1{n-k},\dots,-\frac 1{n-k}
}_{\text{$n-k$ times}}\Big)\,, \quad k=1,\dots,n-1\,.
$$
We estimate the values $l_m(\l^{(k)})$ by comparing sums to 
integrals and obtain
$$
\multline
\left|l_m(\l^{(k)})
\right|\le\frac1k\int_0^k x^m\,dx+\frac1{n-k}\int_k^n x^m\,dx = \\
=\frac1{m+1}\left(k^m+\frac{n^{m+1}-k^{m+1}}{n-k}\right)\le 2n^m
\endmultline 
$$
This concludes the proof. \qed
\enddemo

Now we proceed with the proof of the theorem

\demo{Proof of Theorem 5.4}
It is clear that it suffices to prove the estimate \tht{5.15} for
any set of functions $f_1,f_2,\dots\in\Lt$ whose highest degree
terms generate the algebra $\L$. We take
$$
f_m(\l)=\sum_i \left[(\l_i-\t(i-1))^m - (-\t(i-1))^m\right]\,,
\quad m=1,2,\dots\,.
$$
After we expand the binomials $(\l_i-\t(i-1))^m$ we shall have
to estimate the sums
$$
\left|\sum_{i=1}^n \l_i^r (i-1)^{m-r} \right|\,, \quad r=1,\dots,m\,.
$$
For $r=1$ this was done in Lemma 5.5 which implies
$$
\left|\sum_{i=1}^n \l_i  (i-1)^{m-1} \right| \le n^{m-1}\, \N{\l} +
\frac 2{\sqrt{\t}}\,n^{m-2}\, \N{\l}^2\,.
$$
For $r\ge 2$ we have
$$
\left|\sum_{i=1}^n \l_i^r (i-1)^{m-r} \right| 
\le n^{m-r}\, \sum (\l_i^2)^{r/2} 
\le n^{m-r} \left( \sum \l_i^2\right)^{r/2} \le n^{m-r} \, \N{\l}^r
\,.
$$
This concludes the proof. \qed
\enddemo

\head
6.~Proof of Theorem 1.3
\endhead

In this section we shall give a proof of a general abstract
approximation theorem which will imply that
$$
\Ex \U \subset
\big\{\P_{\a,\b,\g}\big\}_{\a^\pm,\b^\pm,\g^\pm}\,. \tag 6.1 
$$
Similar approximation theorems are well-known in the 
literature (see e.g.\ \cite{V}, \cite{VK1}, and \cite{Ol1}, Section 22). 

The inverse inclusion 
$$
\Ex \U \supset
\big\{\P_{\a,\b,\g}\big\}_{\a^\pm,\b^\pm,\g^\pm}
$$
follows immediately from the fact that the functions 
$\P_{\a,\b,\g}$, being products of the form \tht{1.5}, are
extreme in a larger convex set, namely, in the set of
characteristic functions of $S(\infty)$-invariant measures on
$$
\Z^\infty:=\varprojlim \Z^n \,. \tag 6.2 
$$
Recall that the description of ergodic 
$S(\infty)$-invariant measures on infinite Cartesian products 
like \tht{6.2} is given by the De Finetti theorem, see e.g.\ \cite{Al}. 

Now we formulate the abstract setup of the approximation
theorem. 
Let a set $A$ be presented as a disjoint union 
$$
A=\bigsqcup_{n\ge 1} A_n
$$
of countable sets $A_1,A_2,\dots$. This presentation can
be also encoded as a function
$$
\gather
[\,\cdot\,]:A\to\Z_{>0}\,, \\
[a]=n \Leftrightarrow a\in A_n\,.
\endgather
$$
Suppose also that we are
given  a function
$$
\om:\bigsqcup_{n\ge 1} \left(A_n\times A_{n+1} \right)\to [0,1]
$$
satisfying
$$
\sum_{a\in A_n} \om(a,a')=1\,, \quad \forall n, \quad \forall a'\in
A_{n+1}\,. \tag 6.3
$$
To this data one associates a convex set $\A$ as follows.
Let
$$
\R^{A_n} \cong \R\times\R\times\dots
$$
denote the vector space of real-valued functions
$$
f: A_n \to \R \,.
$$
This is a locally convex metrizable vector space with
respect to the product topology. The subset 
$$
\R^{A_n}\supset \A_n:=\{f:A_n\to[0,1], \sum_{a\in A_n} f(a) = 1 \}
$$
is a simplex whose vertices are the delta-functions
$$
\d_a(a')=\cases 1\,, &a'=a\,,\\
0\,, &a'\ne a\,, \endcases \quad a,a'\in A_n\,.
$$
The function $\om$ defines a projection 
$$
\Proj{n+1}n:\A_{n+1} \to \A_n
$$
by the following formula
$$
\left(\Proj{n+1}n f \right)(a):=\sum_{a'\in A_{n+1}} 
\om(a,a') f(a')\,, \quad a\in A_n, \quad f\in \A_{n+1} \,. \tag 6.4
$$
For any $n>k\ge 1$ we denote by $\Proj nk$ the composition
of the projections \tht{6.4}. Now we set, by definition
$$
\A:=\varprojlim \A_n
$$
with respect to these projections. This is a convex set; any
element of $\A$ is by construction a function on $A$. Denote
by $\Ex\A$ the set of extreme points of $\A$. We shall now
establish the following approximation theorem for functions in  
$\Ex\A$. 

\proclaim{Theorem 6.1} For any function $f\in\Ex\A$ there
exists a sequence 
$$
(a_n)\subset A\,, \quad [a_n]\to\infty\,,
$$
of elements of the set $A$ such that
$$
\lim_{n\to\infty} \Proj{[a_n]}k \d_{a_n} = f\big|_{A_k}\,, \quad \forall
k=1,2,\dots \,.
$$
\endproclaim 

\remark{Remark} In fact, a simple argument shows that we can
always choose the sequence $(a_n)\subset A$ in such a way that
$[a_n]=n$.
\endremark 

\demo{Proof of Theorem 6.1}
We shall follow the argument used in proof of Theorem 22.9 in
\cite{Ol1} and shall consider certain compactifications of the
non-compact (for infinite $A_n$) sets $\A_n$.  Namely, set
$$
\Ab_n:=\{f:A_n\to[0,1], \sum_{a\in A_n} f(a) \le 1 \}
\subset \R^{A_n}\,.
$$
It is clear that $\Ab_n\supset\A_n$ and $\Ab_n$ is convex and
compact. The same formula \tht{6.4} gives  a map
$$
\bProj{n+1}n:\Ab_{n+1}\to\Ab_n\,.
$$
Observe that this map is not, in general, continuous. However, it
is always lower semicontinuous
$$
\bProj{n+1}n(f_0) \le \varliminf_{f\to f_0}  
\bProj{n+1}n(f)\,, \quad f,f_0\in\Ab_{n+1}\,. \tag 6.5
$$
Here and below an inequality between two functions means 
inequality at every point. The semicontinuity \tht{6.5}  
follows from the fact that $\om(a,a')\in[0,1]$ for any $a$ and
$a'$. By virtue of \tht{6.5} the convex set
$$
\{(f,f'), f \ge \bProj{n+1}n f'\} \subset
\Ab_n\times\Ab_{n+1}
$$
is compact. Similarly, the convex sets $\B_n$
$$
\B_n:=\{(f_1,\dots,f_n), f_i \ge \bProj{i+1}i
f_{i+1}, 1\le i<n\} \subset \prod_1^n \Ab_i 
$$
are compact and so are the sets
$$
\fC_n:=\B_n\times\prod_{n+1}^\infty \Ab_i \subset \prod_{1}^\infty
\Ab_i \subset \R^A\,. 
$$
Clearly, we have  
$$
\fC_1 \supset \fC_2 \supset \dots  \,.
$$
Set
$$
\fC:=\bigcap_{n\ge 1}\fC_n \,.
$$
Since all sets $\fC_n$ are compact convex 
subsets of a locally convex metrizable vector space
$\R^A$ of real-valued functions on $A$ 
it is well known (see e.g.\ Lemma 22.13 in \cite{Ol1})
that for any $f\in\Ex\fC$ there exists
a sequence 
$$
(f_n\in\Ex\fC_n)_{n=1,2,\dots}
$$
such that
$$
f_n \to f\,,\quad  n\to \infty \,.
$$ 

One easily obtains (see Sections 22.20--22.23 in \cite{Ol1})
the following properties of the sets  $\Ex\fC_n$ and $\Ex\fC$.
We have:
\roster
\item
for any point $f\in\Ex\fC_n$ its projection onto $\B_n$
(given by restriction to $\cup_{i\le n}A_i$) lies in $\Ex\B_n$\,;
\item 
a point in $\Ex\B_n$ is either zero or of the form
$$
\Delta_a(a'):=\cases 
\left(\Proj{[a]}{[a']}\d_{a}\right)
(a')\,, & [a']\le[a]\,, \\
0\,, & [a']>[a]\,,
\endcases
$$
where $a\in\A$ and $[a]\le n$; 
\item the inclusion $\A\subset\fC$ induces the inclusion
$$
\Ex\A\subset\Ex\fC\,.
$$
\endroster

Now we can finish the proof of the theorem. For any 
$f\in\Ex\A\subset\Ex\fC$ there exists a sequence 
$(f_n\in\Ex\fC_n)_{n=1,2,\dots}$
such that $f_n\to f$. Consider the projections of $f_n$ 
onto $\B_n$. Clearly, only finitely many of them can be zero.
The non-zero ones produce a sequence $(a_n)\subset A$ such that
$$
\lim_{n\to\infty} \Delta_{a_n} = f\,. \tag 6.6
$$
It is clear that \tht{6.6} implies $[a_n]\to\infty$ and that
$(a_n)$ is the desired sequence. This concludes the proof. \qed
\enddemo

Now we shall specialize this general theorem to our particular
situation. We set
$$
A_n:=\{\text{ signatures $\l$ of length $\ell(\l)=n$ }\}\,.
$$
If $\mu$ and $\l$ are two signatures of length $n$ and $n+1$
respectively we set
$$
\multline
{}\qquad \om(\mu,\l):=\text{ coefficient of $\P_\mu(z_1,\dots,z_n;\t)$}\\
\text{in the expansion of $\P_\l(z_1,\dots,z_n,1;\t)$}\,.\qquad {}
\endmultline
$$
The normalization
$$
\P_\nu(1,\dots,1;\t)=1\,, \quad \forall\nu\,,
$$
implies the property \tht{6.3}.
Now by the definition of $\A$ and the definitions of Section 1.3
we have
$$ 
\A_n \cong \Un\,, \quad \text{and}\quad \A\cong \U \,.
$$
Moreover, the convergence in $\A_n$ is equivalent to 
the uniform convergence of functions in $\Un$. Now Theorems 6.1 and
1.1 imply \tht{6.1} and this concludes the proof of Theorem 1.3.

\head 7.~Appendix. A direct proof of the formula \tht{2.10}
for generating functions
\endhead

In this section we give a direct derivation of the formula 
\tht{2.10} for the generating functions 
and show how the formula \tht{2.8} follows
from this computation. Recall that the formula \tht{2.8} reads
$$
g^*_k(x_1,\dots,x_n;\t):=\sum_{1\le i_1 \le \dots \le i_k \le n}
\frac{(\t)_{m_1} \cdots (\t)_{m_n}}{m_1!\, \cdots m_n!} 
(x_{i_1}-k+1) \cdots (x_{i_{k-1}}-1) x_{i_k}\,, \tag 7.1
$$
where $m_l:=\#\{r\,|\, i_r=l\}$ are the multiplicities with
which the numbers $l=1,\dots,n$ occur in $i_1,\dots,i_k$. 
Set  
$$
G^*_n(x_1,\dots,x_n;u):=\sum_{k\ge 0}
\frac{g^*_k(x_1,\dots,x_n)}{u(u-1)\cdots(u-k+1)} \,.
$$
We shall establish the following

\proclaim{Proposition 7.1} For fixed $x_1,\dots,x_n$ the series
converges provided $\Re u \ll 0$ and 
$$
G^*_n(x_1,\dots,x_n;u) = 
\prod_{i=1}^n 
\frac{\G(x_i-u-\t\,i)}{\G(x_i-u-\t\,i+\t)}
\frac{\G(-u-\t\,i)}{\G(-u-\t\,i+\t)}\,. \tag 7.2 
$$
\endproclaim

\demo{Proof} Induct on $n$. For $n=1$
we have the Gauss summation (see \cite{WW})
$$
\sum_{k\ge 0} \frac{(\t)_k}{k!}
\frac{x(x-1)\cdot(x-k+1)}{u(u-1)\cdot(u-k+1)}=
{}_2F_1(-x,\t;-u;1)=
\frac{\G(x-u-\t)}{\G(x-u)}
\frac{\G(-u)}{\G(-u-\t)}\,, \tag 7.3 
$$
provided $\Re u \ll 0$.
Now suppose $n\ge 2$. By \tht{7.1}  we have
$$
g^*_k(x_1,\dots,x_n)=\sum_{p+q=k} g^*_p(x_1-q,\dots,x_{n-1}-q) \,
g^*_q(x_n)\,.
$$
This implies
$$
G^*_n(x_1,\dots,x_n;u)=\sum_{q\ge 0} G^*_{n-1}(x_1-q,\dots,x_{n-1}-q;u-q)
\frac{g^*_q(x_n)}{u(u-1)\cdots(u-q+1)} \,. \tag 7.4 
$$
Since we assume the truth of \tht{7.2} for $G^*_{n-1}$ we conclude that
$$
\align
\frac{G^*_{n-1}(x_1-q,\dots,x_{n-1}-q;u-q)}
{G^*_{n-1}(x_1,\dots,x_{n-1};u)}&=
\prod_{i=1}^{n-1}
\frac{\G(q-u-\t\,i+\t)}{\G(q-u-\t\,i)}
\frac{\G(-u-\t\,i)}{\G(-u-\t\,i+\t)} \\
&=\prod_{i=1}^{n-1} \frac{(-u-\t\,i+\t)_q}{(-u-\t\,i)_q} \\
&=\frac{u\cdots(u-q+1)}{(u+\t(n-1))\cdots(u+\t(n-1)-q+1)}\,. \tag 7.5
\endalign
$$
Substituting \tht{7.5} into \tht{7.4} and 
then using the Gauss summation  we obtain
$$
G^*_n(x_1,\dots,x_n;u)=G^*_{n-1}(x_1,\dots,x_{n-1};u)
\frac{\G(x-u-\t\,n)}{\G(x-u-\t(n-1))}
\frac{\G(-u-\t(n-1))}{\G(-u-\t\,n)}\,,
$$
which proves \tht{7.2}. \qed
\enddemo 

Now we claim that this computation proves, in fact, the formula
\tht{7.1}. Indeed, in the above prove we have used only 
the explicit formula appearing on the RHS of \tht{7.1}. The result
of the computation shows that the polynomial in the RHS of \tht{7.1}
is an element of $\Lt(n)$. One easily checks the degree and
the vanishing conditions and concludes that \tht{7.1} is true.  

The above argument has a $q$-analog (where $q$ is an extra parameter
and not the index used above). It uses Heine's $q$-analog (see \cite{GR})
of the Gauss summation \tht{7.3} which is equivalent to the following
summation (compare to \tht{2.10} in \cite{Ok9})
$$
\sum_{k\ge 0} t^{-k}
\frac{(t;q)_k}{(q;q)_k}
\frac{(x-1)\cdots(x-q^{k-1})}
{(u-1)\cdots(u-q^{k-1})} =
\frac{(x/u\phantom{t};q)_\infty}
{(x/ut;q)_\infty}
\frac{(1/ut;q)_\infty}
{(1/u\phantom{t};q)_\infty}\,.
$$
Here $(a;q)_k=(1-a)\cdots(1-q^{k-1}a)$.

\enddocument

\end